\documentclass[10pt,oneside]{amsart}

\usepackage{amsmath,amssymb,amsfonts,amssymb,amsthm}

\usepackage{verbatim}
\usepackage[usenames]{color}
\usepackage{hyperref,mathabx}
\usepackage{url}
\usepackage{tikz,tikz-qtree,ifthen,cancel}
\usepackage{array,tikz-qtree,ifthen,cancel}
\usepackage{graphicx}
\usepackage{adjustbox}
\usepackage{amsthm,graphicx,tikz,appendix,tikz-qtree,ifthen,cancel}
\usepackage{appendix}
\usepackage{caption, subcaption}

\usetikzlibrary{calc,shapes,patterns,positioning}

\usepackage{yfonts}

\usepackage{soul}

\usepackage{todonotes}

\newtheorem{thm}{Theorem}[section]
\newtheorem{prop}[thm]{Proposition}

\newtheorem*{thmFBRD}{Theorem \ref{thm.FBRD}}
\newtheorem*{thmMain}{Theorem \ref{thmmain}}
\newtheorem*{lemEnd-homog}{Lemma \ref{lem.End-homog}}
\newtheorem{lem}[thm]{Lemma}
\newtheorem{cor}[thm]{Corollary}
\newtheorem{observation}[thm]{Observation}
\newtheorem*{thm7.3}{Theorem \ref{thm.MillikenSWP}}

\theoremstyle{remark}
\newtheorem{rem}[thm]{Remark}

\theoremstyle{definition}
\newtheorem{defn}[thm]{Definition}

\newtheorem{assumption}[thm]{Assumption}

\newtheorem{question}[thm]{Question}

\theoremstyle{remark}

\newcommand{\re}{\upharpoonright}

\newcommand{\al}{\alpha}
\newcommand{\om}{\omega}

\newcommand{\sse}{\subseteq}
\newcommand{\contains}{\supseteq}
\newcommand{\forces}{\Vdash}

\DeclareMathOperator{\NR}{R}

\DeclareMathOperator{\ran}{ran}

\DeclareMathOperator{\height}{ht}
\DeclareMathOperator{\depth}{depth}

\DeclareMathOperator{\Crit}{Crit}
\DeclareMathOperator{\CL}{CL}
\DeclareMathOperator{\IS}{IS}
\DeclareMathOperator{\BRD}{BRD}
\DeclareMathOperator{\EBRD}{EBRD}

\DeclareMathOperator{\Emb}{Emb}

\newcommand{\ra}{\rightarrow}

\newcommand{\Llra}{\Longleftrightarrow}
\newcommand{\lgl}{\langle}
\newcommand{\rgl}{\rangle}

\newcommand{\Erdos}{Erd{\H{o}}s}
\newcommand{\Fraisse}{Fra{\"{i}}ss{\'{e}}}

\newcommand{\Lauchli}{L{\"{a}}uchli}
\newcommand{\Wazewski}{Wa\.{z}ewski}

\newcommand{\bP}{\mathbb{P}}

\newcommand{\bS}{\mathbb{S}}

\newcommand{\perfect}{ray}

\newcommand{\noprint}[1]{\relax}

\title[Big Ramsey degrees and the two-branching pseudotree]{Big Ramsey degrees and the two-branching pseudotree}

\author{David Chodounsk\'{y}}
\address{Institute of Mathematics of the Czech Academy of Sciences,
Žitná~25, Praha~1, Czech Republic, and
Department of Applied Mathematics (KAM), Charles University, Malostranské náměstí~25, Praha~1,
Czech Republic.
}
\email{chodounsky@math.cas.cz}

\author{Natasha Dobrinen}
\address{Department of Mathematics\\
 University of Notre Dame \\
255 Hurley Bldg\\
Notre Dame, IN 46556   U.S.A.}
\email{ndobrine@nd.edu}

\author{Thilo Weinert}
\address{
University of Vienna \\
Institute for Mathematics\\
Kurt G\"odel Research Centre for Mathematical Logic\\
Kolingasse 14--16 \\
1090 Wien\\
Austria}
\email{thilo.weinert@univie.ac.at}

\thanks{
Chodounsk\'{y} was supported by project 25-15571S of the Czech Science Foundation (GAČR) and by the Czech Academy of Sciences CAS (RVO 67985840). 
Dobrinen was supported by National Science Foundation Grant DMS-2300896.
Weinert was supported by the Fonds zur F\"orderung der wissenschaftlichen Forschung, Lise Meitner grant M 3037-N, as well as the Research Project of National relevance ``PRIN2022\_DIMONTE - Models, sets and classifications (realizzato con il contributo del progetto PRIN 2022 - D.D. n. 104 del 02/02/2022 – PRIN2022\_DIMONTE - Models, sets and classifications - Codice 2022TECZJA\_003 - CUP N. G53D23001890006. “Finanziato dall'Unione Europea – Next-GenerationEU – M4 C2 I1.1”). For the purpose of open access, the authors have applied a CC BY public copyright licence to any Author Accepted Manuscript version arising from this submission.
}


\begin{document}

\maketitle
\tableofcontents

\begin{abstract}
We prove that  each finite chain in the two-branching countable ultrahomogeneous pseudotree has finite big Ramsey degrees. This is in contrast to the recent result of Chodounsk\'{y}, Eskew, and Weinert that antichains of size two have infinite big Ramsey degree in  the pseudotree.  Combining a lower bound result of theirs with work in this paper shows that chains of length two in the pseudotree have big Ramsey degree exactly seven.  The  pseudotree  is the first example of a countable ultrahomogeneous structure in a finite language in which some finite substructures have finite big Ramsey degrees while others have  infinite big Ramsey degrees. 
\end{abstract}

\section{Introduction}

This paper 
contributes to   a rapidly expanding line of research on big Ramsey degrees of infinite structures
and 
presents the first positive
results on an ultrahomogeneous structure whose age is a  \Fraisse\ $-$ HP class with a Ramsey expansion. 
We prove that finite chains in the  two-branching countable ultrahomogeneous pseudotree, which we denote by  $\Psi$, have finite big Ramsey degrees. 
This is in contrast to the recent  result in \cite{CEW_EUROCOMB25} that antichains of size at least two in  $\Psi$  have infinite big Ramsey degrees.
Thus, the homogeneous two-branching pseudotree  has  finite big Ramsey degrees for certain finite substructures and infinite big Ramsey degrees for others. 
This is the first example of an ultrahomogeneous structure in a finite language with this behavior. 

To introduce the area of big Ramsey degrees, we begin with  Ramsey's celebrated result.

\begin{thm}[Infinite Ramsey Theorem, \cite{Ramsey30}]\label{thm.RamseyInfinite}
Given positive integers $m$ and $r$,
suppose the  collection of  all $m$-element subsets of $\mathbb{N}$
is colored by  $r$ colors.
Then there is an infinite set  $N$ of natural numbers
such that all $m$-element subsets of $N$  have the same color.
\end{thm}


When moving to infinite structures, exact analogues of Ramsey's theorem are often not possible.
This was first seen in 1933 in the paper \cite{Sierpinski} of Sierpi\'{n}ski, in which he demonstrated a two-coloring of pairs of rationals such that every subcopy of the rationals retains both colors. 
In unpublished work, Galvin showed that for each finite coloring of pairs of rationals, there is a subcopy of the rationals and two colors so that each pair in the subcopy takes one of those two colors.
This was the earliest observation of what is now known as big Ramsey degrees.

\begin{defn}[\cite{Kechris/Pestov/Todorcevic05}]\label{defn.bRd}
Given an infinite structure $\mathbf{S}$ and a finite substructure $\mathbf{A}$ of $ \mathbf{S}$,
let $\BRD(\mathbf{A},\mathbf{S})$ denote the least integer $n\ge 1$, if it exists, such that
given any coloring of 
the copies of $\mathbf{A}$ in $\mathbf{S}$ into finitely many colors, 
 there is a
 substructure $\mathbf{S}'$ of $\mathbf{S}$, isomorphic to $\mathbf{S}$,  such that 
the copies of $\mathbf{A}$ in 
$\mathbf{S}'$ take no more than $n$ colors.
Otherwise, we write $\BRD(\mathbf{A},\mathbf{S})=\infty$.
We say that
 $\mathbf{S}$ has {\em finite big Ramsey degrees} if for each finite substructure $\mathbf{A}$ of $\mathbf{S}$,
$\BRD(\mathbf{A},\mathbf{S})\in\mathbb{N}$.
\end{defn}

In \cite{Kechris/Pestov/Todorcevic05},  
Kechris, Pestov, and Todorcevic proved 
a  correspondence between the Ramsey property of a \Fraisse\ class and the extreme amenability of the automorphism group of its \Fraisse\ structure.
Additionally, they demonstrated how  big Ramsey degrees for \Fraisse\  structures $\mathbf{F}$ are related to
big oscillation degrees for their automorphism groups, Aut$(\mathbf{F})$, and asked for more examples of homogeneous structures with finite big Ramsey degrees.

Historically, big Ramsey degrees have been defined in terms of substructures but more recently, especially in connection with topological dynamics, 
big Ramsey degrees have also been defined  in terms of colorings of  embeddings.
The {\em embedding big Ramsey degree} of $\mathbf{A}$ in $ \mathbf{S}$, denoted here by  $\EBRD(\mathbf{A},\mathbf{S})$,
is  the least integer $n\ge 1$, if it exists, such that
given any coloring of 
the embeddings  of $\mathbf{A}$ in $\mathbf{S}$ into finitely many colors, 
 there is a
 substructure $\mathbf{S}'$ of $\mathbf{S}$, isomorphic to $\mathbf{S}$,  such that 
the embeddings of $\mathbf{A}$ in 
$\mathbf{S}'$ take no more than $n$ colors.
There is a simple connection between the two:
$\EBRD(\mathbf{A},\mathbf{S}) =|\mathrm{Aut}(\mathbf{A})|\cdot\BRD(\mathbf{A},\mathbf{S})$,
where $\mathrm{Aut}(\mathbf{A})$ is the set of automorphisms of $\mathbf{A}$.
A theorem of Hjorth in 
\cite{Hjorth08} showed that whenever a homogeneous structure has more than one automorphism,  there is  some finite substructure of it with embedding  big Ramsey degree at least two (possibly infinite).
Zucker proved  in \cite{Zucker19} that
if  a  \Fraisse\ structure  $\mathbf{F}$ has finite big Ramsey degrees and moreover, $\mathbf{F}$  admits a big Ramsey structure,
then
any big Ramsey flow of Aut$(\mathbf{F})$ is a universal completion flow, and further,  any two universal completion flows are isomorphic.

Extrapolating on specific questions asked in \cite{Kechris/Pestov/Todorcevic05} leads to the following natural question (mentioned as Question 5.1 in \cite{DobrinenIfCoLog20}):
If  a \Fraisse\ class has  a precompact Ramsey expansion, does its \Fraisse\ limit have finite big Ramsey degrees?
Sauer already
 negatively answered this question 
in \cite{Sauer03} by exhibiting a  
 \Fraisse\ free amalgamation class of directed graphs, with countably infinitely many forbidden substructures,
whose \Fraisse\ limit  has infinite big Ramsey degree for vertices.
Recently,
 Braunfeld, Chodounsk\'{y},  de Rancourt,  Hubi\v{c}ka,  Kawach, and Kone\v{c}n\'{y} 
showed  in \cite{BCdRHKK24}
that the answer is negative 
for all \Fraisse\ free amalgamation classes 
where  the signature is allowed  infinitely many relations with the same arity $n\ge 2$.
Very recently,
work of Chodounsk\'{y}, Eskew, and Weinert in \cite{CEW_EUROCOMB25} showed  that in the pseudotree,
singletons have big Ramsey degree one, whereas
antichains of size two  or more  have infinite big Ramsey degree;  this is in a structure whose signature has only finitely many relations.
In contrast, we show in this paper that all finite chains in the pseudotree have finite big Ramsey degrees.

The first  infinite structures in which big Ramsey degrees were fully characterized are 
the rationals (D.~Devlin in \cite{DevlinThesis})
and the Rado graph (moreover, all binary relational simple structures) 
 via the two papers (Sauer \cite{Sauer06} and Laflamme, Sauer, and Vuksanovic \cite{Laflamme/Sauer/Vuksanovic06}).
Other work  between 2006 and 2010 includes
\cite{NVT08},  \cite{NVTSauer09}, and \cite{Laflamme/NVT/Sauer10}.
Work on big Ramsey degrees halted after this due to lack of methods viable for structures with forbidden substructures. 
Then in 2017,  the second author completed a proof that   the triangle-free Henson graph has finite big Ramsey degrees in \cite{DobrinenJML20}, and that all $k$-clique-free Henson graphs have finite big Ramsey degrees in \cite{DobrinenJML23} (completed in 2018) 
answering  questions in \cite{Kechris/Pestov/Todorcevic05}.
These works 
 developed new methods for dealing with homogeneous structures which have  forbidden substructures, 
  instigating  a spate of research on big Ramsey degrees, currently totalling over  45 papers  by a multitude of authors  since  then.
For background and references, 
see the expository  paper \cite{Dobrinen_ICM}  
connected with the second author's 2022 ICM invited lecture, and the recent survey paper \cite{HZ25}.

Motivation for our current investigation of big Ramsey degrees of the pseudotree comes from work of Kwiatkowska on the Ramsey theory and topological dynamics of dendrites  in \cite{Kwiatkowska18}, presented at the 2018 Banff workshop on {\em Unifying Themes in Ramsey Theory}.
A {\em continuum} is a compact  connected topological space.
A {\em dendrite} is a locally connected continuum that contains no simple closed curve.
The {\em order} of a point $x$ is the number of connected components obtained after removing $x$.  A {\em ramification point} is a point of order $\ge 3$. 
For $P\sse\{3,4,\dots,\om\}$,
a generalized  {\em \Wazewski\ dendrite} $W_P$ is a dendrite such that each ramification point has order in $P$, and each arc  in $W_P$ contains a ramification point.  

For every $P\sse\{3,4,\dots,\om\}$, 
Charatonik and  Dilks proved in \cite{CharatonikDilks94} that 
$W_P$ is unique up to homeomorphism. 
Duchesne and Monod 
proved in \cite{DuchesneMonod19} that 
the homeomorphism groups of generalized \Wazewski\ dendrites are simple.
Kwiatkowska  investigated the universal minimal flow of the  homeomorphism group of $W_P$.

\begin{thm}[Kwiatkowska, \cite{Kwiatkowska18}]\label{thm.KUMF}
If $P\sse \{3,4,\dots,\om\}$ is finite, then the universal minimal flow of the homeomorphism group $H(W_P)$ is metrizable, and is computed explicitly.
If $P$ is infinite, then the universal minimal flow  of $H(W_P)$ is not metrizable.
\end{thm}

Let  $M_P$ denote   the set of all ramification points of $W_P$.  Then  $M_P$ is  a countable structure which is dense in $W_P$. 
Duchesne showed in 
 \cite{Duchesne20} that 
$M_P$ is ultrahomogeneous with respect to 
its age.
Duchesne and Monod showed in \cite{DuchesneMonod19} that 
H$(W_P)$,  with the uniform metric, is isomorphic (as a top.\ group) to  Aut$(M_P)$, with the pointwise convergence metric.  
For  recent work on projective \Fraisse\ limits and dendrites, see 
\cite{CodenottiKwiatkowska24}.
 For more related work, see the paper \cite{BJP16} by Bodirsky, Jonsson, Pham  and 
 the paper \cite{KRS21} by
Kaplan, Rzepecki, and Siniora.

Let $\mathcal{F}_P$ denote the  set of all finite  substructures of $M_P$.
The members of  $\mathcal{F}_P$  are finite unordered trees where the degree of each vertex is in $\{1\}\cup P$.
Kwiatkowska proved in Proposition 1.2 in \cite{Kwiatkowska18} that $M_P$ is ultrahomogeneous   with respect to  $\mathcal{F}_P$.
She then
proved and utilised the following theorem in order  to prove Theorem \ref{thm.KUMF}.

\begin{thm}[Kwiatkowska, \cite{Kwiatkowska18}]\label{thm.KRE}
$\mathcal{F}_P$ has a precompact Ramsey expansion.
\end{thm}

Theorem \ref{thm.KRE}  naturally leads to the following question, which was asked at the 2018 Banff Workshop on {\em Unifying Themes in Ramsey Theory}.

\begin{question}\label{q.1.5}
For each finite $P\sse \{3,4,\dots,\om\}$, does $M_P$ have finite big Ramsey degrees?
\end{question}

The structures  $M_P$   are closely connected with pseudotrees.
Designating one node in $M_P$ as the root  induces a partial ordering $\preceq$ on $M_P$.  
Upon deleting the root, each remaining connected component  of $M_P$, with the induced partial order $\preceq$, is a pseudotree. 
Seen another way,
intermediate between $\mathcal{F}_P$ and Kwiatkowska's Ramsey expansion 
 $\mathcal{F}^*_P$  is the class of 
finite trees $\mathcal{T}_P$ with nodes of branching degree in $\{1\}\cup P$ and  tree partial order $\preceq$.
For structures  $A,B\in \mathcal{T}_P$,   we say that $A$ embeds into $B$ iff there is a $\preceq$-preserving one-to-one map from the nodes in $A$ into $B$. Then the \Fraisse\ limit of $\mathcal{T}_P$ is a pseudotree, which we denote by $\Psi_P$.  This leads to the natural modification of Question \ref{q.1.5}.

\begin{question}
For each finite $P\sse \{3,4,\dots,\om\}$, does $\Psi_P$ have finite big Ramsey degrees?
\end{question}

To begin this investigation, we have focused on  the special case $P=\{3\}$,  letting 
 $\Psi$ denote $\Psi_{\{3\}}$, 
the two-branching countable pseudotree (see Definition \ref{def.Psi}).
As mentioned above, Chodounsk\'{y}, Eskew, and Weinert recently proved the following:

\begin{thm}[Chodounsk\'{y}, Eskew, and Weinert, \cite{CEW_EUROCOMB25}]\label{thm.CEW}
Letting $\mathbf{A}$ denote a singleton,
$\BRD(\mathbf{A},\Psi)=1$.
On the other hand, 
if $\mathbf{A}$ is an antichain in $\Psi$ of size two or more,
then $\BRD(\mathbf{A},\Psi)=\infty$.
\end{thm}

In contrast, we show the following:

\begin{thm}[Main Theorem]\label{thmmain}
For each finite chain $\mathbf{C}$ in $\Psi$, 
$\BRD(\mathbf{C},\Psi)<\infty$.
In particular, $\BRD(\mathbf{C},\Psi)$ is no more than the number of diaries for $\mathbf{C}$.
\end{thm}

Chodounsk\'{y}, Eskew, and Weinert further proved in \cite{CEW_EUROCOMB25} that chains of length two have big Ramsey degree at least seven.
In Corollary \ref{cor.BRD7}, we show that the big Ramsey degree for chains of length two is at most seven; hence, the degree is exactly seven.

The paper is organized as follows.
In Section \ref{sec.basics} we define our particular version of  coding trees for the pseudotree  $\Psi$.
Included is a discussion of  the intuition behind these coding trees,
their merits for aiding us in proving Ramsey theorems for $\Psi$, and several observations translating between $\Psi$ and coding trees.
Section \ref{sec.diaries}
introduces diaries and 
 provides an outline of the proof of the Main Theorem.
 Section \ref{sec.HLVariant} builds the
Ramsey-theoretic tools on coding trees to be applied later in the proof of  the Main Theorem.  
This includes amalgamation lemmas and
new variants of the Halpern--\Lauchli\ and Milliken Theorems for diaries in coding trees for $\Psi$.  
The
Main Theorem is proved in Section \ref{sec.mainresults}.
This is via first finding a minimalistic subtree of a coding tree which still encodes a copy of $\Psi$, and then applying the theorems from Section  \ref{sec.HLVariant} to achieve finite upper bounds.
We conclude that section by exhibiting that there are only seven 
diaries for chains of length two
two inside $\Psi$, showing that their big Ramsey degree is at most seven.


\section{Coding trees for $\Psi$}\label{sec.basics}

This section provides background setting up the structures on which we will prove Ramsey theorems leading to upper bounds for big Ramsey degrees of finite chains in $\Psi$.
In particular, it contains 
the expanded structure $\Psi^*$ (which is intrinsic in its big Ramsey degrees),
the definition of  our coding trees,  and  observations to help the reader interpret between   $\Psi^*$ and coding trees.

 Pseudotrees have appeared in  the literature using different languages, (e.g., \cite{Kurepa77}, \cite{MonkBK}, \cite{KRS21}).
 For ease of connecting the results in this paper with those in \cite{CEW_EUROCOMB25}, we repeat  their definition of pseudotree here.

\begin{defn}\label{def.Psi}
A  {\em pseudotree} is a structure $\Psi$ 
in language $\mathcal{L}=\{\preceq, \wedge\}$
where 
$\preceq$ 
is an order (reflexive, antisymmetric, transitive) 
and 
 $\wedge$ is  a binary function satisfying 
$p\wedge q=\sup\{t\in \Psi:t\preceq p,\ t\preceq q\}$ for each $p,q\in \Psi$ (the supremum always exists). 
Moreover, for each $p\in \Psi$, the set $D(p)=\{t\in \Psi:t\prec p\}$ is linearly ordered by $\prec$, where  $t\prec p$ if and only if  $t\preceq p$ and $t\ne p$.
\end{defn}

The pseudotree under investigation in this paper is the 
{\em binary} or {\em $2$-branching universal pseudotree}, which we will denote by
$\Psi$.
Up to isomorphism,  
$\Psi$ is the countable pseudotree with the properties that
 for each $p\in\Psi$, $D(p)$ is order-isomorphic to the rationals and 
there are distinct $q,r\in\Psi\setminus\{p\}$ with $p=q\wedge r$, while
for  any three nodes $s,t,u\in \Psi$, it is not the case that $p=s\wedge t=s\wedge u=t\wedge u$.
Equivalently,  $\Psi$  is the \Fraisse\ limit of the class of finite trees   $\mathcal{T}$ in the language $\mathcal{L}=\{\preceq, \wedge\}$ where each non-terminal node in  $T$ has at most two immediate successors.

Towards finding upper bounds for (and  exactly characterizing in a forthcoming paper) the big Ramsey degrees of $\Psi$, we 
 expand the language to $\mathcal{L}^*=\mathcal{L}\cup\{<_{\mathrm{lex}}\}$ 
and define $\mathcal{T}^*$ to be the class of finite trees
$T$ in the language $\mathcal{L}^*$ 
such that 
the reduct of $T$ to $\mathcal{L}$ is in $\mathcal{T}$, and 
$<_{\mathrm{lex}}$ is a linear order on $T$
satisfying the following:
The order $<_{\mathrm{lex}}$ orders each pair of 
 $\prec$-immediate successors of some $t\in T$.
Further,
if  $p$ and $q$ are  $\preceq$-incomparable in $T$ and $t=p\wedge q$, then 
 $p<_{\mathrm{lex}} q$ if and only if $p'<_{\mathrm{lex}} q'$ where $p',q'$ are the immediate successors of $t$ so that  $p'\preceq p$ and $q'\preceq q$; in this case, 
$p<_{\mathrm{lex}} t$ and $t<_{\mathrm{lex}} q$.
This order provides  a notion of left and right extensions in $T$.
Note that the leftmost branch in $T$ is the maximal $\prec$-chain $C \sse T$ such that 
for any two distinct nodes $p,q\in C$, 
$p\prec q$ if and only if $q<_{\mathrm{lex}}p$.

Let  $\Psi^*=(P;\preceq,\wedge, <_{\mathrm{lex}})$ denote pseudotree $\Psi$ with the additional order $<_{\mathrm{lex}}$.
Note that the age of $\Psi^*$ equals $\mathcal{T}^*$, but $\Psi^*$ is not the \Fraisse\ limit of $\mathcal{T}^*$.
For each node $p\in P$, 
the sets 
\begin{align}
\Psi^*(p,0)&:=\{q\in P: p\prec q\mathrm{\ and\ } q<_{\mathrm{lex}}p\}\cr
 \Psi^*(p,1)&:=\{q\in P:p\prec q\mathrm{\ and\ } p<_{\mathrm{lex}}q\}
\end{align}
 partition the nodes in $P$ that are $\prec$-above $p$ into two copies of $\Psi^*$, the `left' and `right' pseudotrees above $p$, respectively.
Let 
\begin{equation}
\Psi^*(p,2):=\{q\in P: \neg(p\prec q)\}.
\end{equation}
Then $\Psi^*(p,2)$ is exactly the set of those $q$ whose meet with $p$ is $\prec$-below $p$.
Notice that the sets
$\Psi^*(p,0)$, $\Psi^*(p,1)$, $\Psi^*(p,2)$ partition 
$ \Psi^*\setminus\{p\}$ into three copies of $\Psi^*$.
Given  $p\in P$,  let $L_p$ denote
the {\em leftmost ray in $\Psi^*$ above $p$}, defined to be the
maximal $\prec$-chain $L_p\sse 
 \Psi^*(p,0)$ such that  
for any $q,r\in  L_p$,
$q\prec r$   implies $r<_{\mathrm{lex}}q$.
Notice that $L_p$ is a dense linear order without endpoints.

All big Ramsey degree results  involve an enumeration of the infinite structure in order type $\om$.  This essential feature  was  discovered by Sierpi\'{n}ski in  \cite{Sierpinski}, which showed that there is a coloring of pairs of rationals into two colors, each of which persists in every subcopy of the rationals. 
The enumeration of a structure induces a tree of $1$-types, a standard notion in model theory. 
Coding trees  were first formulated in 
\cite{DobrinenJML20} 
 and have  been  used in 
\cite{DobrinenJML23}, 
\cite{DobrinenRado19}, \cite{Zucker22},
 \cite{Balko7},  \cite{CDP1}, \cite{CDP2},  \cite{Dobrinen_SDAP}, and \cite{Dobrinen/Zucker23} to prove big Ramsey degree and infinite-dimensional Ramsey theorems on a large collection of \Fraisse\ structures. 
The coding trees   in 
all of the above-mentioned papers are trees
of  $1$-types that are realizable and complete (over finite initial segments of the enumerated structure).
In this paper, we use a slightly different construction.  Our coding trees will have nodes that correspond to  $1$-types that are realizable but not always complete.
In fact, using the tree of complete $1$-types would be unnecessarily cumbersome and would not have the following property, which is useful in our proofs:
In our coding tree (Definition \ref{defn.codingtree1} below), each level consists of nodes that 
partition $\Psi^*$ (minus finitely many nodes) into finitely may pieces, each isomorphic to $\Psi^*$.
Such partitions are determined
via  intersections of finitely many
sets of the form $\Psi^*(p,k)$, for $k\in\{0,1,2\}$;
these correspond to sets of all nodes realizing some $1$-type (not necessarily complete) with finitely many parameters.
We will  discuss  further benefits of and motivations behind this set-up throughout this section.

We mention here that big Ramsey degrees are, by definition,  independent of any enumeration of the structure.
An enumeration is an inherent factor in diaries, and further,  in big Ramsey structures (see \cite{Zucker19} and \cite{HZ25}).
For a given homogeneous structure $\mathbf{K}$, the extension property  allows one to construct, within one enumerated copy of $\mathbf{K}$ a subcopy of $\mathbf{K}$ with any other given enumeration.
This has been well-known for  some time, and was made known to the second author via conversations with Sauer in 2014.
See also  Theorem 4.1 of \cite{Masulovic18}, where Ma\v{s}ulovi\'{c} showed, for all countably infinite relational structures  (not necessarily homogeneous) with ages that have strong amalgamation, finiteness of  big Ramsey degrees does not depend on the enumeration.
See also   Corollary 3.17 of \cite{HZ25} for \Fraisse\ classes with strong amalgamation.
For classes without HP, it makes the most sense to enumerate the structure in such a way that all finite initial segments of the enumeration are algebraically closed.

There are several reasons for   using coding trees to obtain our results.
First, an enumeration of the vertices in an infinite structure forms an essential component  of  characterizations of big Ramsey degrees,
so we may as well start with one.  
Second, the coding tree structure makes possible Ramsey theorems in the vein of the  Halpern--\Lauchli\ and Milliken
 theorems.
Third, the proofs of exact degrees and  of infinite-dimensional Ramsey theorems necessarily, in one form or another, involve coding trees.
Our presentation in this paper sets the stage for that future work. 
(See \cite{Dobrinen_ICM} for a discussion of the different methods for proving big Ramsey degrees and their various merits and drawbacks.)

We now define a particular coding tree $\bS$ which will provide the template for our space of coding trees.
We will want to enumerate the nodes in $\Psi^*$ in order-type $\om$ in such a way that every initial segment of this order  is a  member of $\mathcal{T}^*$.
For $S,T\in\mathcal{T}^*$, $S$ {\em embeds} into $T$ iff there is an injection $f:S\ra T$ preserving $\prec$, $\wedge$, and $\le_{\mathrm{lex}}$.
Such an $f$ is an {\em embedding} of $S$ into $T$.
We say that $S$ is a {\em substructure} of $T$ iff  the identity map on $S$ is an embedding of $S$ into $T$.
Although  any enumeration for which each  finite intial segment  forms a tree will suffice to achieve our results,
for concreteness we fix a particular enumeration of $\Psi^*$ which induces a coding tree with a regular structure. 
We hope this makes the ideas more quickly accessible to the reader.

\begin{defn}[Coding Tree $\bS$ for $\Psi^*$]\label{defn.codingtree1}
We fix a sequence   $\langle T_n:n<\om\rangle$
of  members of $\mathcal{T}^*$
as in Figures 1.\ and  3.,  where 
 for each $n<\om$, 
$|T_n|=n+1$, 
 $T_n$ is a substructure of $T_{n+1}$ 
and $\bigcup_{n<\om}T_n=\Psi^*$.
For  $n<\om$,  $p_{n+1}$ denotes the unique  node in $T_{n+1}$ that is not in $T_n$.
This produces an enumeration $\langle p_n:n<\om\rangle$ of the nodes in $\Psi^*$ where  for each $n<\om$, the substructure of $\Psi^*$ restricted to the nodes in $\{p_i:i<n+1\}$  is exactly the tree $T_n$.

Our coding tree $\bS$ will be a subtree of $\{0,1,2\}^{<\om}$ with distinguished coding node $c_n$ of length $n$ representing $p_n$.
For any node $s\in \bS$ we let $|s|$ denote the {\em length} of $s$.
We simplify notation by  letting $(p_n,k)$ denote 
 $\Psi^*(p_n,k)$, for $k\in\{0,1,2\}$.
Each $p_n$ 
  introduces a partition of $\Psi^*\setminus\{p_n\}$ into the three pieces $(p_{n},k)$, $k\in\{0,1,2\}$, each of which is isomorphic to $\Psi^*$.
This will correspond to each coding node $c_n$ in  $\bS$ having three immediate successors.
The non-coding nodes in $\bS$ will have a single immediate successor, because they do not introduce any partition of $\Psi^*$.
We now define $\bS$ recursively.  See Figures  1.--4.\ for intuition.

Let 
$c_0$  be the empty sequence, representing $p_0$.  
In  our coding  tree $\bS\sse \{0,1,2\}^{<\om}$,
$c_0$ has three immediate successors $\langle 0\rangle, \langle 1\rangle,  \langle 2\rangle, $  representing the sets $(p_0,0)$,  $(p_0,1)$,  and $(p_0,2)$, respectively.
That is, for each $n\ge 1$, $p_n\in (p_0,k)$ if and only if the coding node $c_n$ extends $\langle k\rangle$.

Let $c_1$ denote the node $\langle 0\rangle$; this 
corresponds to the fact that $p_1\in (p_0,0)$, as $p_0\prec p_1$ and $p_1<_{\mathrm{lex}}p_0$ in $T_1$.
The immediate successors of $c_1$ are ${c_1}^{\frown}\langle k\rangle$, $k\in\{0,1,2\}$.
The node $\langle 0,0\rangle$ represents  $(p_0,0)\cap(p_1,0)$,  
the set of all $p_n$, $n>1$, satisfying $p_1\prec p_n$ and $p_n<_{\mathrm{lex}}p_1$.
That is, every node  $p_n\in (p_0,0)\cap(p_1,0)$ will be represented by a  coding node  $c_n$ in $\bS$ end-extending $\langle 0,0\rangle$.
Similarly,
$\langle 0,1\rangle$ represents 
 $(p_0,0)\cap(p_1,1)$, 
the set of all $p_n$, $n>1$, satisfying $p_1\prec p_n$ and $p_1<_{\mathrm{lex}}p_n$;
and 
$\langle 0,2\rangle$ represents 
 $(p_0,0)\cap(p_1,2)$, 
the set of all $p_n$, $n>1$,
 satisfying $p_0\prec p_n\prec p_1$.
All $p_n\in (p_0,1)\cup (p_0,2)$ are automatically in $(p_1,2)$.
Thus, the immediate successors of $\langle 1\rangle$ and $\langle 2\rangle$ in $\bS$ are necessarily $\langle 1,2\rangle$ and $\langle 2,2\rangle$, respectively.
We point out that the nodes in $\bS\cap \{0,1,2\}^2$ 
correspond to the five pieces that 
the set of nodes $\{p_0,p_1\}$ partitions $\Psi^*\setminus \{p_0,p_1\}$ into.
Let $c_2=\langle 1,2\rangle$, since $p_2\in (p_0,1)\cap(p_1,2)$. 
This concludes the construction of the nodes in $\bS$ up to level $2$.

The following conditions are preserved in the recursive construction of $\bS$.
\begin{enumerate}
\item
The  node $p_n\in\Psi^*$ is uniquely represented by the coding node $c_n$ of length $n$, where 
 $c_n$  denotes the sequence  $\langle k_i:i<n\rangle$  in $\{0,1,2\}^n$   such that  for each $i<n$, 
$p_n\in (p_i,k_i)$.
\item
The splitting nodes in $\bS$ are exactly the coding nodes in $\bS$, and each has three immediate successors.
\item
Every non-coding node  has a unique extension, either by $0$ or by $2$.
\item
For each $n\ge 1$, the $n$-th level of $\bS$ has exactly $1 +2^n$ many nodes, representing the  partition of $\Psi^*$ into $1+2^n$ many pieces induced by $\{p_0,\dots,p_{n-1}\}$.
That is, $s=\langle \ell_i:i<n\rangle\in \bS$ represents $\bigcap_{i<n}(p_i,\ell_i)$.
\end{enumerate}

For $n\ge 2$,  suppose we have  the tree $\bS$ up to level $n$ satisfying (1)--(4).
Then 
 $c_n$  denotes the sequence  $\langle k_i:i<n\rangle$  in $\{0,1,2\}^n$   such that  for each $i<n$, 
$p_n\in (p_i,k_i)$.
The node $p_n$ partitions $\bigcap_{i<n}(p_i,k_i)$ into three new pieces, namely
$(\bigcap_{i<n}(p_i,k_i))\cap (p_n,k)$ for $k\in\{0,1,2\}$,
 so the node $c_n$ has three immediate successors, $c_n^{\frown}\langle k\rangle$, $k\in\{0,1,2\}$, in $\bS$.

For   $s=\langle \ell_i:i<n\rangle\in \bS$, a non-coding node at level $n$,
$s$ represents 
$\bigcap_{i<n}(p_i,\ell_i)$.
Note that $p_n$ does not partition this piece, so $s$ has a unique extension in $\bS$.
To see this in detail, 
let $j<n$ be such that $s\wedge c_n=c_j$.
 We have three cases:

\underline{Case 1}. $c_n$ extends ${c_j}^{\frown}0$.
 Then  $p_j\prec p_n$, and 
$s$ extends ${c_j}^{\frown}1$ or ${c_j}^{\frown}2$.
Either way, any $p_m$ in $\bigcap_{i<n}(p_i,\ell_i)$
satisfies $p_m\wedge p_n\preceq p_j\prec p_n$, so $p_m\in (p_n,2)$. 
Hence, 
 $s$  must be extended by $s^{\frown}2$.

\underline{Case 2}.
$c_n$ extends ${c_j}^{\frown}1$.
Then  $p_j\prec p_n$.
If $s$ extends ${c_j}^{\frown}0$, then 
any $p_m$  in $\bigcap_{i<n}(p_i,\ell_i)$
 satisfies  $p_m\wedge p_n=p_j$. 
If $s$ extends ${c_j}^{\frown}2$, then 
any $p_m$  in $\bigcap_{i<n}(p_i,\ell_i)$
 satisfies  $p_m\wedge p_n=p_j$. 
Either way, $p_m\in (p_n,2)$, so $s$ must be extended by $s^{\frown}2$.

\underline{Case 3}.  $c_n$ extends ${c_j}^{\frown}2$. 
Then  $p_n\wedge p_j\prec p_j$, 
and $s$ extends ${c_j}^{\frown}0$ or ${c_j}^{\frown}1$ so any 
 $p_m$  in $\bigcap_{i<n}(p_i,\ell_i)$ satisfies $p_j\prec p_m$. 
If $p_n\prec p_j$, then  also $p_n\prec p_m$ so $p_m\in (p_n,0)$ and $s$ is extended by $s^{\frown}0$.
Note that in this case, $c(i)\in\{0,2\}$ for each  $j\le i<n$.
Otherwise, 
$p_n\not\prec p_j$ so $p_m\wedge p_n\prec p_n$ and $s$ is extended by $s^{\frown}2$.
Note that in this case, there is some $j<i<n$ for which $c_n(i)=1$. 

Let $c_{n+1}$ denote the node $\langle k'_i:i<n+1\rangle$  in level $n+1$  of $\bS$ such that 
$p_{n+1}$ is in $\bigcap_{i<n+1}(p_i,k'_i)$.
This concludes the construction of level $n+1$ of $\bS$, satisfying (1)--(4).
\end{defn}

 Figure 1.\  exhibits the  first four trees in our sequence $\langle T_n:n<\om\rangle$.
Figure 2.\ shows  the first four levels of the coding tree $\bS$ induced by these four trees.
Figure 3.\ shows  our enumeration of the first $13$ nodes in $\Psi^*$, in particular, $T_{13}$ in our sequence, and Figure 4.\ shows the coding tree induced by this enumeration.

\begin{figure}
\begin{subfigure}[b]{0.2\textwidth}
\centering
\resizebox{\linewidth}{!}{
\begin{tikzpicture}[scale=.15]
\foreach \x in {0}{
\foreach \y in {0}{
\node  at (\x,\y) {};
}}
\foreach \x in {0,4,8,12}{
\foreach \y in {-8,0,4,6,8,12,16}{
\node at (\x,\y) {};
}}
\node[circle, fill=black,inner sep=0pt, minimum size=5pt] at (0,0) {};
\node[left] at (0,0) {$p_0$};
\end{tikzpicture}
}
\caption{$T_0$}
\label{fig:subfigA}
\end{subfigure}
\begin{subfigure}[b]{0.2\textwidth}
\centering
\resizebox{\linewidth}{!}{
\begin{tikzpicture}[scale=.15]
\foreach \x in {0}{
\foreach \y in {0}{
\node  at (\x,\y) {};
}}
\foreach \x in {0,4,8,12}{
\foreach \y in {-8,0,4,6,8,12,16}{
\node at (\x,\y) {};
}}
\draw (0,0)--(0,8);
\node[circle, fill=black,inner sep=0pt, minimum size=5pt] at (0,0) {};
\node[left] at (0,0) {$p_0$};
\node[circle, fill=black,inner sep=0pt, minimum size=5pt] at (0,8) {};
\node[left] at (0,8) {$p_1$};
\end{tikzpicture}
}
\caption{$T_1$}
\label{fig:subfigB}
\end{subfigure}
\begin{subfigure}[b]{0.2\textwidth}
\centering
\resizebox{\linewidth}{!}{
\begin{tikzpicture}[scale=.15]
\foreach \x in {0}{
\foreach \y in {0}{
\node  at (\x,\y) {};
}}
\foreach \x in {0,4,8,12}{
\foreach \y in {-8,0,4,6,8,12,16}{
\node at (\x,\y) {};
}}
\draw (0,0)--(0,8);
\draw(0,0)--(8,8);
\node[circle, fill=black,inner sep=0pt, minimum size=5pt] at (0,0) {};
\node[left] at (0,0) {$p_0$};
\node[circle, fill=black,inner sep=0pt, minimum size=5pt] at (0,8) {};
\node[left] at (0,8) {$p_1$};
\node[circle, fill=black,inner sep=0pt, minimum size=5pt] at (8,8) {};
\node[left] at (8,8) {$p_2$};
\end{tikzpicture}
}
\caption{$T_2$}
\label{fig:subfigC}
\end{subfigure}
\begin{subfigure}[b]{0.2\textwidth}
\centering
\resizebox{\linewidth}{!}{
\begin{tikzpicture}[scale=.15]
\foreach \x in {0}{
\foreach \y in {0}{
\node  at (\x,\y) {};
}}
\foreach \x in {0,4,8,12}{
\foreach \y in {-8,0,4,6,8,12,16}{
\node at (\x,\y) {};
}}
\draw (0,-8)--(0,0)--(0,8);
\draw(0,0)--(8,8);
\node[circle, fill=black,inner sep=0pt, minimum size=5pt] at (0,-8) {};
\node[left] at (0,-8) {$p_3$};
\node[circle, fill=black,inner sep=0pt, minimum size=5pt] at (0,0) {};
\node[left] at (0,0) {$p_0$};
\node[circle, fill=black,inner sep=0pt, minimum size=5pt] at (0,8) {};
\node[left] at (0,8) {$p_1$};
\node[circle, fill=black,inner sep=0pt, minimum size=5pt] at (8,8) {};
\node[left] at (8,8) {$p_2$};
\end{tikzpicture}
}
\caption{$T_3$}
\label{fig:subfigD}
\end{subfigure}
\caption{A  sequence of finite trees in $\mathcal{T}^*$}
\end{figure}


\begin{figure}
\begin{tikzpicture}[scale=.25]

\foreach \x in {0}{
\foreach \y in {0}{
\node  at (\x,\y) {};
}}

\foreach \x in {-18,0,18}{
\foreach \y in {2}{
\node at (\x,\y) {};
}}
\draw (-18,2)--(0,0)--(18,2);
\draw(0,0)--(0,2);

\node[circle, fill=black,inner sep=0pt, minimum size=5pt] at (0,0) {};
\node[left] at (0.7,-1)  {$c_0$};

\foreach \x in {-24,-18,-12,0,18}{
\foreach \y in {4}{
\node at (\x,\y) {};
}}
\draw (-24,4)--(-18,2)--(-12,4);
\draw(-18,2)--(-18,4);
\draw(0,2)--(0,4);
\draw(18,2)--(18,4);
\node[circle, fill=black,inner sep=0pt, minimum size=5pt] at (-18,2) {};
\node[left] at (-18,1.5) {$c_1$};
\node[left] at (-5.5,0.2){$(p_0,0)$};
\node[left] at (2.2,1){$(p_0,1)$};
\node[left] at (8.6,0.15){$(p_0,2)$};
\node[left] at (-19.6,2.5){$(p_1,0)$};
\node[left] at (-15.7,4.2){$(p_1,1)$};
\node[left] at (-11.3,2.5){$(p_1,2)$};

\foreach \x in {-24,-18,-12,-6,0,6,18}{
\foreach \y in {6}{
\node at (\x,\y) {};
}}
\draw (-24,4)--(-24,6);
\draw(-18,4)--(-18,6);
\draw(-12,4)--(-12,6);

\draw(0,6)--(0,4);
\draw(-6,6)--(0,4)--(6,6);
\draw(18,4)--(18,6);

\node[circle, fill=black,inner sep=0pt, minimum size=5pt] at (0,4) {};
\node[left] at (0,3.8) {$c_2$};
\node[left] at (-2.1,4.5){$(p_2,0)$};
\node[left] at (2.2,6){$(p_2,1)$};
\node[left] at (6.5,4.5){$(p_2,2)$};

\foreach \x in {-24,-18,-12,-6,0,6,12,18,24}{
\foreach \y in {8}{
\node at (\x,\y) {};
}}
\draw (-24,8)--(-24,6);
\draw(-18,8)--(-18,6);
\draw(-12,8)--(-12,6);

\draw(-6,8)--(-6,6);
\draw(0,6)--(0,8);
\draw(6,8)--(6,6);

\draw(18,8)--(18,6);
\draw(12,8)--(18,6)--(24,8);

\node[circle, fill=black,inner sep=0pt, minimum size=5pt] at (18,6) {};
\node[left] at (18,5.8) {$c_3$};
\node[left] at (16,6.5){$(p_3,0)$};
\node[left] at (20.3,7.4){$(p_3,1)$};
\node[left] at (24.8,6.5){$(p_3,2)$};

\end{tikzpicture}
\caption{First $4$ levels of $\mathbb{S}$ with labeled $1$-types}
\end{figure}


\begin{figure}
\begin{tikzpicture}[scale=.2]

\foreach \x in {0}{
\foreach \y in {0}{
\node  at (\x,\y) {};
}}

\foreach \x in {0,4,8,16}{
\foreach \y in {-16,-8, -4,0,4,6,8,12,16}{
\node at (\x,\y) {};
}}
\draw (0,-8)--(0,0)--(0,8);
\draw(0,0)--(8,8);

\node[circle, fill=black,inner sep=0pt, minimum size=5pt] at (0,-8) {};
\node[left] at (0,-8) {$p_3$};
\node[circle, fill=black,inner sep=0pt, minimum size=5pt] at (0,0) {};
\node[left] at (0,0) {$p_0$};
\node[circle, fill=black,inner sep=0pt, minimum size=5pt] at (0,8) {};
\node[left] at (0,8) {$p_1$};
\node[circle, fill=black,inner sep=0pt, minimum size=5pt] at (8,8) {};
\node[left] at (8,8) {$p_2$};

\draw (0,16)--(0,8);
\draw(0,8)--(8,16);
\draw(16,16)--(8,8);
\draw(16,12)--(8,8);

\node[circle, fill=black,inner sep=0pt, minimum size=5pt] at (0,16) {};
\node[left] at (0,16) {$p_4$};
\node[circle, fill=black,inner sep=0pt, minimum size=5pt] at (8,16) {};
\node[left] at (8,16) {$p_5$};
\node[circle, fill=black,inner sep=0pt, minimum size=5pt] at (0,4) {};
\node[left] at (0,4) {$p_6$};
\node[circle, fill=black,inner sep=0pt, minimum size=5pt] at (16,16) {};
\node[left] at (16,16) {$p_7$};
\node[circle, fill=black,inner sep=0pt, minimum size=5pt] at (16,12) {};
\node[left] at (17.2,10.8) {$p_8$};
\node[circle, fill=black,inner sep=0pt, minimum size=5pt] at (4,4) {};
\node[left] at (4,4) {$p_9$};

\node[circle, fill=black,inner sep=0pt, minimum size=5pt] at (0,-4) {};
\node[left] at (0,-4) {$p_{10}$};

\draw (0,-8)--(8,0);
\node[circle, fill=black,inner sep=0pt, minimum size=5pt] at (8,0) {};
\node[left] at (8,0) {$p_{11}$};

\draw (0,-16)--(0,-8);
\node[circle, fill=black,inner sep=0pt, minimum size=5pt] at (0,-16) {};
\node[left] at (0,-16) {$p_{12}$};

\end{tikzpicture}
\caption{Enumerating the first 13 nodes in $\Psi$, or $T_{13}$}
\end{figure}


\begin{figure}
\begin{tikzpicture}[scale=.21]

\foreach \x in {0}{
\foreach \y in {0}{
\node  at (\x,\y) {};
}}

\foreach \x in {-18,0,18}{
\foreach \y in {2}{
\node at (\x,\y) {};
}}
\draw (-18,2)--(0,0)--(18,2);
\draw(0,0)--(0,2);

\node[circle, fill=black,inner sep=0pt, minimum size=5pt] at (0,0) {};
\node[left] at (0.7,-1)  {$c_0$};

\foreach \x in {-24,-18,-12,0,18}{
\foreach \y in {4}{
\node at (\x,\y) {};
}}
\draw (-24,4)--(-18,2)--(-12,4);
\draw(-18,2)--(-18,4);
\draw(0,2)--(0,4);
\draw(18,2)--(18,4);

\node[circle, fill=black,inner sep=0pt, minimum size=5pt] at (-18,2) {};
\node[left] at (-18,1.5) {$c_1$};

\foreach \x in {-24,-18,-12,-6,0,6,18}{
\foreach \y in {6}{
\node at (\x,\y) {};
}}
\draw (-24,4)--(-24,6);
\draw(-18,4)--(-18,6);
\draw(-12,4)--(-12,6);

\draw(0,6)--(0,4);
\draw(-6,6)--(0,4)--(6,6);
\draw(18,4)--(18,6);

\node[circle, fill=black,inner sep=0pt, minimum size=5pt] at (0,4) {};
\node[left] at (0,3.8) {$c_2$};

\foreach \x in {-24,-18,-12,-6,0,6,12,18,24}{
\foreach \y in {8}{
\node at (\x,\y) {};
}}
\draw (-24,8)--(-24,6);
\draw(-18,8)--(-18,6);
\draw(-12,8)--(-12,6);

\draw(-6,8)--(-6,6);
\draw(0,6)--(0,8);
\draw(6,8)--(6,6);

\draw(18,8)--(18,6);
\draw(12,8)--(18,6)--(24,8);

\node[circle, fill=black,inner sep=0pt, minimum size=5pt] at (18,6) {};
\node[left] at (18,5.8) {$c_3$};

\foreach \x in {-26,-24,-22,-18,-12,-6,0,6,12,18,24}{
\foreach \y in {10}{
\node at (\x,\y) {};
}}
\draw (-24,8)--(-24,10);
\draw(-26,10)--(-24,8)--(-22,10);

\draw(-18,8)--(-18,10);
\draw(-12,8)--(-12,10);

\draw(-6,8)--(-6,10);
\draw(0,10)--(0,8);
\draw(6,8)--(6,10);

\draw(12,8)--(12,10);
\draw(18,8)--(18,10);
\draw(24,8)--(24,10);

\node[circle, fill=black,inner sep=0pt, minimum size=5pt] at (-24,8) {};
\node[left] at (-24,7.8) {$c_4$};

\foreach \x in {-26,-24,-22,-20,-18,-16,-12,-6,0,6,12,18,24}{
\foreach \y in {12}{
\node at (\x,\y) {};
}}
\draw(-26,10)--(-26,12);
\draw (-24,12)--(-24,10);
\draw(-22,12)--(-22,10);

\draw(-18,12)--(-18,10);
\draw(-20,12)--(-18,10)--(-16,12);

\draw(-12,12)--(-12,10);

\draw(-6,12)--(-6,10);
\draw(0,10)--(0,12);
\draw(6,12)--(6,10);

\draw(12,12)--(12,10);
\draw(18,12)--(18,10);
\draw(24,12)--(24,10);

\node[circle, fill=black,inner sep=0pt, minimum size=5pt] at (-18,10) {};
\node[left] at (-18,9.8) {$c_5$};

\foreach \x in {-26,-24,-22,-20,-18,-16,-14,-12,-10,-6,0,6,12,18,24}{
\foreach \y in {14}{
\node at (\x,\y) {};
}}
\draw(-26,14)--(-26,12);
\draw (-24,12)--(-24,14);
\draw(-22,12)--(-22,14);

\draw(-16,12)--(-16,14);
\draw(-18,12)--(-18,14);
\draw(-20,12)--(-20,14);

\draw(-12,12)--(-12,14);
\draw(-14,14)--(-12,12)--(-10,14);

\draw(-6,12)--(-6,14);
\draw(0,14)--(0,12);
\draw(6,12)--(6,14);

\draw(12,12)--(12,14);
\draw(18,12)--(18,14);
\draw(24,12)--(24,14);

\node[circle, fill=black,inner sep=0pt, minimum size=5pt] at (-12,12) {};
\node[left] at (-12,11.8) {$c_6$};

\foreach \x in {-26,-24,-22,-20,-18,-16,-14,-12,-10,-8,-6,-4,0,6,12,18,24}{
\foreach \y in {16}{
\node at (\x,\y) {};
}}
\draw(-26,14)--(-26,16);
\draw (-24,16)--(-24,14);
\draw(-22,16)--(-22,14);

\draw(-20,16)--(-20,14);
\draw(-18,16)--(-18,14);
\draw(-16,16)--(-16,14);

\draw(-14,14)--(-14,16);
\draw(-12,14)--(-12,16);
\draw(-10,14)--(-10,16);

\draw(-6,16)--(-6,14);
\draw(-8,16)--(-6,14)--(-4,16);

\draw(0,14)--(0,16);
\draw(6,16)--(6,14);

\draw(12,16)--(12,14);
\draw(18,16)--(18,14);
\draw(24,16)--(24,14);

\node[circle, fill=black,inner sep=0pt, minimum size=5pt] at (-6,14) {};
\node[left] at (-6,13.8) {$c_7$};

\foreach \x in {-26,-24,-22,-20,-18,-16,-14,-12,-10,-8,-6,-4,-2,0,2,6,12,18,24}{
\foreach \y in {18}{
\node at (\x,\y) {};
}}
\draw(-26,18)--(-26,16);
\draw (-24,16)--(-24,18);
\draw(-22,16)--(-22,18);

\draw(-20,16)--(-20,18);
\draw(-18,16)--(-18,18);
\draw(-16,16)--(-16,18);

\draw(-14,18)--(-14,16);
\draw(-12,18)--(-12,16);
\draw(-10,18)--(-10,16);

\draw(-8,16)--(-8,18);
\draw(-6,16)--(-6,18);
\draw(-4,16)--(-4,18);

\draw(0,18)--(0,16);
\draw(-2,18)--(0,16)--(2,18);

\draw(6,16)--(6,18);

\draw(12,16)--(12,18);
\draw(18,16)--(18,18);
\draw(24,16)--(24,18);

\node[circle, fill=black,inner sep=0pt, minimum size=5pt] at (0,16) {};
\node[left] at (0,15.8) {$c_8$};

\foreach \x in {-26,-24,-22,-20,-18,-16,-14,-12,-10,-8,-6,-4,-2,0,2,4,6,8,12,18,24}{
\foreach \y in {20}{
\node at (\x,\y) {};
}}
\draw(-26,18)--(-26,20);
\draw (-24,20)--(-24,18);
\draw(-22,20)--(-22,18);

\draw(-20,20)--(-20,18);
\draw(-18,20)--(-18,18);
\draw(-16,20)--(-16,18);

\draw(-14,18)--(-14,20);
\draw(-12,18)--(-12,20);
\draw(-10,18)--(-10,20);

\draw(-8,20)--(-8,18);
\draw(-6,20)--(-6,18);
\draw(-4,20)--(-4,18);

\draw(-2,18)--(-2,20);
\draw(0,18)--(0,20);
\draw(2,18)--(2,20);

\draw(6,20)--(6,18);
\draw(4,20)--(6,18)--(8,20);

\draw(12,20)--(12,18);
\draw(18,20)--(18,18);
\draw(24,20)--(24,18);

\node[circle, fill=black,inner sep=0pt, minimum size=5pt] at (6,18) {};
\node[left] at (6,17.8) {$c_9$};

\foreach \x in {-26,-24,-22,-20,-18,-16,-14,-12,-10,-8,-6,-4,-2,0,2,4,6,8,10,12,14,18,24}{
\foreach \y in {22}{
\node at (\x,\y) {};
}}
\draw(-26,22)--(-26,20);
\draw (-24,20)--(-24,22);
\draw(-22,20)--(-22,22);

\draw(-20,20)--(-20,22);
\draw(-18,20)--(-18,22);
\draw(-16,20)--(-16,22);

\draw(-14,22)--(-14,20);
\draw(-12,22)--(-12,20);
\draw(-10,22)--(-10,20);

\draw(-8,20)--(-8,22);
\draw(-6,20)--(-6,22);
\draw(-4,20)--(-4,22);

\draw(-2,22)--(-2,20);
\draw(0,22)--(0,20);
\draw(2,22)--(2,20);

\draw(4,20)--(4,22);
\draw(6,20)--(6,22);
\draw(8,20)--(8,22);

\draw(12,20)--(12,22);
\draw(10,22)--(12,20)--(14,22);

\draw(18,20)--(18,22);
\draw(24,20)--(24,22);

\node[circle, fill=black,inner sep=0pt, minimum size=5pt] at (12,20) {};
\node[left] at (12,19.8) {$c_{10}$};

\foreach \x in {-26,-24,-22,-20,-18,-16,-14,-12,-10,-8,-6,-4,-2,0,2,4,6,8,10,12,14,16,18,20,24}{
\foreach \y in {24}{
\node at (\x,\y) {};
}}
\draw(-26,22)--(-26,24);
\draw (-24,20)--(-24,24);
\draw(-22,20)--(-22,24);

\draw(-20,20)--(-20,24);
\draw(-18,20)--(-18,24);
\draw(-16,20)--(-16,24);

\draw(-14,22)--(-14,24);
\draw(-12,22)--(-12,24);
\draw(-10,22)--(-10,24);

\draw(-8,24)--(-8,22);
\draw(-6,24)--(-6,22);
\draw(-4,24)--(-4,22);

\draw(-2,22)--(-2,24);
\draw(0,22)--(0,24);
\draw(2,22)--(2,24);

\draw(4,24)--(4,22);
\draw(6,24)--(6,22);
\draw(8,24)--(8,22);

\draw(10,22)--(10,24);
\draw(12,22)--(12,24);
\draw(14,22)--(14,24);

\draw(18,24)--(18,22);
\draw(16,24)--(18,22)--(20,24);

\draw(24,24)--(24,22);

\node[circle, fill=black,inner sep=0pt, minimum size=5pt] at (18,22) {};
\node[left] at (18,21.8) {$c_{11}$};

\foreach \x in {-26,-24,-22,-20,-18,-16,-14,-12,-10,-8,-6,-4,-2,0,2,4,6,8,10,12,14,16,18,20,22,24,26}{
\foreach \y in {26}{
\node at (\x,\y) {};
}}
\draw(-26,26)--(-26,24);
\draw (-24,26)--(-24,24);
\draw(-22,26)--(-22,24);

\draw(-20,26)--(-20,24);
\draw(-18,26)--(-18,24);
\draw(-16,26)--(-16,24);

\draw(-14,26)--(-14,24);
\draw(-12,26)--(-12,24);
\draw(-10,26)--(-10,24);

\draw(-8,24)--(-8,26);
\draw(-6,24)--(-6,26);
\draw(-4,24)--(-4,26);

\draw(-2,26)--(-2,24);
\draw(0,26)--(0,24);
\draw(2,26)--(2,24);

\draw(4,24)--(4,26);
\draw(6,24)--(6,26);
\draw(8,24)--(8,26);

\draw(10,26)--(10,24);
\draw(12,26)--(12,24);
\draw(14,26)--(14,24);

\draw(16,24)--(16,26);
\draw(18,24)--(18,26);
\draw(20,24)--(20,26);

\draw(24,24)--(24,26);
\draw(22,26)--(24,24)--(26,26);

\node[circle, fill=black,inner sep=0pt, minimum size=5pt] at (24,24) {};
\node[left] at (24,23.8) {$c_{12}$};

\end{tikzpicture}
\caption{The first 13 levels of $\mathbb{S}$}
\end{figure}


\begin{rem}\label{rem.partition}
For any node $s\in\bS$,  there is an isomorphic copy (see Definition \ref{defn.ctiso}) of $\bS$ in the subtree
of $\bS$ above $s$.
This is crucial to the proofs of Ramsey theorems  in Section \ref{sec.HLVariant}.
\end{rem}

\begin{rem}\label{rem.notCT1T}
Recall  that $(p,0)$ and $(p,1)$ are sets of nodes realizing  a 
 complete $1$-type over $\Psi^*\re \{p\}$, while
$(p,2)$ is the union of the finitely many such sets.
Hence,  each node $s\in \bS$ corresponds to
 a realizable $1$-type over $\Psi^*\re\{p_i:i<|s|\}$, but  not complete if $s(\ell)=2$ for some $\ell<|s|$.
Further, since the $(p,k)$, $k\in\{0,1,2\}$, partition $\Psi^*\setminus \{p\}$ into three copies of $\Psi^*$,
we see that the nodes of length $n$ in $\bS$ partition
$\Psi^*\setminus\{p_i:i<n\}$ into $1+2^n$ many pieces each of which is a copy of $\Psi^*$.
For example, looking at Figures 1.\ and 3.,
the nodes in level $2$ of $\bS$ are $\langle 0,0\rangle$, 
$\langle 0,1\rangle$, $\langle 0,2\rangle$, $\langle 1,2\rangle$, and $\langle 2,2\rangle$.
These nodes correspond to  the five pieces   that $T_1$ partitions 
 $\Psi^*\setminus \{p_0,p_1\}$ into:
\begin{enumerate}
\item
 $\langle 0,0\rangle$ corresponds to $(p_0,0)\cap (p_1,0)$,  the set of  those nodes in $\Psi^*$ that are $\prec$-above and $<_{\mathrm{lex}}$-below   $p_1$  (so above and left of $p_1$).
In this case, the information regarding $p_0$ is redundant.
\item
$\langle 0,1\rangle$ 
corresponds to $(p_0,0)\cap (p_1,1)$,   the set of  those nodes $q$  in $\Psi^*$ that are $\prec$-above  $p_1$ and  $<_{\mathrm{lex}}$-above  $p_1$ (so above and right of $p_1$).
Again, the information regarding $p_0$ is redundant.
\item
$\langle 0,2\rangle$  corresponds to $(p_0,0)\cap (p_1,2)$,  the set of   those nodes $q$  in $\Psi^*$ that are $\prec$-above  $p_0$  with  $q\wedge p_1\prec p_1$ (so those nodes that are either $\prec$ between $p_0$ and $p_1$ or that branch off right  from this region).
\item
 $\langle 1,2\rangle$
corresponds to 
$(p_0,1)\cap (p_1,2)$,  the set of   those nodes $q$  in $\Psi^*$ that are $\prec$-above  $p_0$
with $p_0<_{\mathrm{lex}}q$ (so those nodes that are above and right of $p_0$).
In this case, any $q \in (p_0,1)$
 necessarily satisfies  
 $q\wedge p_1=p_0\prec p_1$, so   the immediate successor of  $\langle 1\rangle$ in $\bS$ must be $\langle 1,2\rangle$.
\item
 $\langle 2,2\rangle$
corresponds to 
$(p_0,2)\cap (p_1,2)$,  the set of those nodes $q$  in $\Psi^*$ that for which $q\wedge p_0\prec p_0$.
Again, any $q\in (p_0,2)$ 
 necessarily satisfies  
 $q\wedge p_1=p_0\prec p_1$, so   the immediate successor of $\langle 2\rangle$ in $\bS$ must be $\langle 2,2\rangle$.
\end{enumerate}
\end{rem}

To  keep track of the (leftmost) rays in $\Psi$ in the order that they appear in   the generic sequence $\langle T_n:n<\om\rangle$, 
we define a function $\theta:\Psi\ra\om$ as follows:
The {\em $0$-th   ray} is  $\NR_0=\{p_n: p_n\preceq p_0\}\cup L_{p_0}$.  Define $\theta(p_n)=0$ for each $p_n\in \NR_0$.
Thus, in Figure 3., $p_0, p_1,p_3,p_4, p_6,p_{10}$ all have $\theta$-value $0$.
The  next node in our enumeration that is not in $\NR_0$ is $p_2$.
The  {\em $1$-st 
 ray}
is $\NR_1=\{p_n: p_0\prec p_n\preceq p_2\}\cup L_{p_2}$. 
Define $\theta(p_n)=1$ for each $p_n\in \NR_1$.
In  Figure 3.,  $p_2,p_7,p_9$ all have $\theta$-value $1$.
In general, for $\ell\ge 2$, 
having defined $\NR_i$ for all $i<\ell$, 
take  $k$ least such that $p_k\not \in \bigcup_{i<\ell} \NR_i$, and
 let $j,m<k$ be such that $p_k\wedge  p_m=p_j$ and $p_m$ is the $\prec$-immediate successor of $p_j$ in $T_{k-1}$.
(Such $j$ and $m$ always exist and are unique by virtue of properties of the sequence $\langle T_n:n<\om\rangle$.)
Then define $\NR_{\ell}=\{p_n: p_j\prec p_n\preceq p_k\}\cup L_{p_k}$, and 
let $\theta(p_n)=\ell$ for each $p_n\in \NR_{\ell}$.

Extend $\theta$ to a function on all nodes of $\bS$ as follows:
For a coding node $c_n\in\bS$, define $\theta(c_n)=\theta(p_n)$.
For a non-coding node $s\in \bS$, let $n$ be least such that the coding node $c_n$ extends $s$ in $\bS$ and define $\theta(s)=\theta(c_n)$.
It is straightforward to check  that for any $n$,   every $t$  that extends $c_n$ by only $0$'s and $2$'s has the same $\theta$-value as $c_n$, leading to the following fact.

\begin{observation}\label{obs.key}
$c_k$ is the coding node of minimal length extending some ${c_i}^{\frown}1$ iff $\NR_{\theta(c_k)} = \{p_n: c_n={c_k}^{\frown}s,\ s\in\{0,2\}^{<\om}\}$.
\end{observation}

\begin{rem}\label{rem.Devlin}
Observation \ref{obs.key} was one of the motivations 
behind using the sets  $(p,k)$ to build  $\bS$.
Given any $s\in \mathbb{S}$, the subtree of $\bS$ consisting of all extensions of $s$ by sequences in $\{0,2\}^{<\om}$ corresponds to a copy of the rationals within some ray in $\Psi$.
For each $n$  for which  $\theta(p_n)\not \in \{\theta(p_0),\dots,\theta(p_{n-1})\}$,
the  ray $\NR_{\theta(p_n)}$  consists exactly of those $p_i$, $i\ge n$, for which $c_i$ extends $c_n$ by only $0$'s and $2$'s.
The subtree of $\bS$ consisting of these $c_i$ is isomorphic to a coding  tree of complete $1$-types for  the rationals.
This is especially useful in  relating diaries  for substructures of $\Psi$ with the well-known Devlin-types, that is,  diaries for subsets of the rationals.
\end{rem}


The next observation  is straightforward from the definitions and included for the convenience of the reader.

\begin{observation}[Translating  from  relations in $\bS$ to $\Psi^*$]\label{obs.easy}
Fix $i\ne j$.
\begin{enumerate}
\item[$\bullet$]
If  ${c_i}^{\frown}0\sse c_j$, then  $p_i\prec p_j$ and $p_j<_{\mathrm{lex}} p_i$.
If further $c_j={c_i}^{\frown}0^{\frown}s$
 where  $s\in\{0,2\}^{<\om}$, then
$p_i$ and $p_j$ are on the same ray.
Otherwise,   they are on different rays.

\item[$\bullet$]
If ${c_i}^{\frown}1\sse c_j$, then  $p_i\prec p_j$, $p_i<_{\mathrm{lex}} p_j$, and  they are on different rays.

\item[$\bullet$]
If 
 $c_j= {c_i}^{\frown}2^{\frown}s$ where $s\in \{0,2\}^{<\om}$, then  $p_j\prec p_i$, $p_i<_{\mathrm{lex}} p_j$, and they are on the same ray.
\item[$\bullet$]
If 
 $c_j= {c_i}^{\frown}2^{\frown}s$  where $s\not\in \{0,2\}^{<\om}$, then
 $p_i$ and  $p_j$ are $\prec$-incompatible and 
  $p_i<_{\mathrm{lex}} p_j$.
\end{enumerate}
For the next points, assume that $c_i$ and $c_j$ are incompatible in $\bS$ and let $k$ be such that $c_k=c_i\cap c_j$.
\begin{enumerate}

\item[$\bullet$]
If ${c_k}^{\frown}0\sse c_i$ and ${c_k}^{\frown}1\sse c_j$, then $p_i$ and $p_j$ are $\prec$-incompatible, $p_i\wedge p_j=p_k$, and 
$p_i<_{\mathrm{lex}} p_j$.

\item[$\bullet$]
If ${c_k}^{\frown}1\sse c_i$ and $c_j= {c_k}^{\frown}2^{\frown}s$ where $s\in \{0,2\}^{<\om}$, then  $p_j\prec p_i$, $p_i<_{\mathrm{lex}}p_j$, and   they are on different rays.

\item[$\bullet$]
If ${c_k}^{\frown}1\sse c_i$ and $c_j= {c_k}^{\frown}2^{\frown}s$ where $s\not\in \{0,2\}^{<\om}$,
then $p_i$ and $p_j$ are $\prec$-incompatible and 
$p_i<_{\mathrm{lex}} p_j$.

\item[$\bullet$]
If ${c_k}^{\frown}0\sse c_i$ and $c_j={c_k}^{\frown}2^{\frown}s$ where $s\in \{0,2\}^{<\om}$, then
$p_j\prec p_k\prec p_i$, $p_i<_{\mathrm{lex}}p_k<_{\mathrm{lex}}p_j$.
Further,  $p_i$ and $p_j$ are in the same ray iff
$c_i=  {c_k}^{\frown}0^{\frown}t$ where $t\in \{0,2\}^{<\om}$.

\item[$\bullet$]
If ${c_k}^{\frown}0\sse c_i$  and $c_j={c_k}^{\frown}2^{\frown}s$  where $s\not\in \{0,2\}^{<\om}$, then
$p_i$ and $p_j$ are $\prec$-incompatible and 
$p_i<_{\mathrm{lex}} p_j$.

\end{enumerate}
\end{observation}

We now describe  $\bS$ as a first-order structure, with use of $\om$.
The language $\mathcal{L}_{CT}=\{\sse, \cap, +, b,\theta, |\cdot|\}$ for our coding trees
has  a binary relation symbol  $\sse$,  ternary relation symbols $\cap, +$, and unary function symbols  $b, \theta, |\cdot|$ which are interpreted in $\bS$, with $(\om,<)$,  as follows:
$\sse$ is the tree relation
 that $\bS$ inherits as a subtree of  $3^{<\om}$.
The function
$b:\bS\ra \{1,2,3\}$ outputs the branching number $b(s)$ in $\bS$,
 and 
$\theta:\mathbb{S}\ra \om$ is the function defined above.
For $s\in \bS$,
$|s|$ is the length of $s$.
The relation  $\cap(s,t,u)$ holds iff $u$ is the meet of $s$ and $t$;  when  $\cap(s,t,u)$ holds
we usually just write $s\cap t=u$.
The ternary relation $+\sse \bS\times\bS\times\{0,1,2\}$ is defined by 
$+(t,s,k)$ holds iff 
 $|s|<|t|$ and $t(|s|)=k$.
This relation yields what has historically been called the {\em passing number} of $t$ at $s$.

We point out  that the $n$-th coding node $c_n$ in $\bS$ is the $n$-th node  $s$  of  the set of nodes in $\bS$  for which $b(s)=3$, 
ordered by increasing  length $|s|$.
Hence, coding nodes are definable in the structure $\langle \bS, \om; \sse, \cap, +, b,\theta, |\cdot|, < \rangle$, where $<$ denotes the usual order on $\om$.  The standard lexicographic order on $\bS$ is definable in this structure via $\cap$ and $+$.

\begin{defn}[Coding Tree Isomorphism]\label{defn.ctiso}
Let $S$ and $T$  be  substructures of $\bS$. 
A function $f:S\ra T$ is a {\em coding tree isomorphism} iff  $f$ is a bijection and for all $s,t,u\in S$, the following hold:
\begin{enumerate}
    \item[(a)]
$s\sse t\longleftrightarrow f(s)\sse f(t)$;
\item[(b)]
$s\cap t=u \longleftrightarrow f(s)\cap f(t)=f(u)$;
\item[(c)]
$+(s,t,k)$ holds in $\bS$ $\longleftrightarrow$  $+(f(s),f(t),k)$ holds in $S$;
\item[(d)]
$b(s)=b(f(s))$;
\item[(e)]
$\theta(s)\le \theta(t)  \longleftrightarrow \theta(f(s))\le \theta(f(t))$;
\item[(f)]
$|s|\le |t|  \longleftrightarrow |f(s)|\le |f(t)|$.
\end{enumerate}
When there is a coding tree isomorphism between $S$ and $T$, we simply say that $S$ and $T$ are  {\em isomorphic}.
\end{defn}

Note that if  $f:\bS\ra S$ 
 is a coding tree isomorphism,
then (a)--(d) along with (f) in the definition above imply that  $f$ is a tree isomorphism in the usual sense, preserving relative lengths of nodes, lexicographic order as subtrees of $3^{<\om}$, and branching numbers.
This implies that 
$f$ takes the $n$-th coding node  $c_n$ in $\bS$ to the $n$-th coding node  (ordered by relative lengths) in $S$, since coding nodes are exactly those nodes with branching number three.

\begin{observation}\label{obs.2.9}
Let $S$ be an isomorphic substructure of $\bS$ and let $\langle i_n:n<\om\rangle$ be the increasing enumeration such that $\{c_{i_n}:n<\om\}$ is the set of coding nodes in $S$.
Then $\Psi^*\re \{p_{i_n}:n<\om\}$ is a subcopy of $\Psi^*$, that is, an infinite homogeneous two-branching pseudotree. 
\end{observation}

\begin{proof}

The proof is a straightforward application of Definitions \ref{defn.codingtree1} and 
\ref{defn.ctiso}.
Let $f:\bS\ra S$ denote the coding tree isomorphism from $\bS$ to $S$.
Recall that for $m<\om$, $\Psi^*\re\{p_n:n\le m\}=T_m$.
Let $U_m$ denote $\Psi^*\re\{p_{i_n}:n\le m\}$.
We will show by induction on $m<\om$ that $U_m$ is isomorphic to $T_m$, as structures in $\mathcal{T}^*$.   The observation then follows immediately.

Before beginning the induction proof,  we point out the following:
For each $n<\om$,   $f(c_n)=c_{i_n}$.
By (e) in  Definition \ref{defn.ctiso}, $\theta(c_n)=\theta(c_j)$ iff $\theta(c_{i_n})=\theta(c_{i_j})$.
By (c) and (f) in Definition \ref{defn.ctiso}, it follows that 
$f$ preserves passing numbers; hence, for $j<n$, $c_n(|c_j|)=c_{i_n}(|c_{i_j}|)$.
Recall  from Definition \ref{defn.codingtree1} that for $j<n$, $c_n(|c_j|)=0$ iff 
$p_j\prec p_n$ and $p_n<_{\mathrm{lex}} p_j$ in $\Psi^*$;
$c_n(|c_j|)=1$ iff 
$p_j\prec p_n$ and $p_j<_{\mathrm{lex}} p_n$ in $\Psi^*$;
and $c_n(|c_j|)=2$ iff 
$\neg(p_j\prec p_n)$  in $\Psi^*$.
Moreover,   if $c_n(|c_j|)=2$, then 
$p_n\prec p_j$ iff $\theta(c_n)=\theta(c_j)$.

The base case $m=0$ is clear, since $T_0=\{p_0\}$ and $U_0=\{p_{i_0}\}$ are both singleton structures in $\mathcal{T}^*$.
Now suppose $m\ge 1$ and $T_{m-1}$  and   $U_{m-1}$ are  isomorphic as structures in $\mathcal{T}^*$, witnessed by  the  map $p_j\mapsto p_{i_j}$ for $j<m$.
We have four cases for how $p_m$ may be related to the nodes in $T_{m-1}$.

\underline{Case 1}. 
$p_m$ is a $\prec$-maximal node  in $T_m$ and  $p_m<_{\mathrm{lex}} p_j$, where 
 $p_j$ is the immediate $\prec$-predecessor of $p_m$ in $T_m$.
That is, $p_m$ is a maximal node in $T_m$ that is on a ray already existing in $T_{m-1}$.
Then  $c_m(|c_j|)=0$ and 
$\theta(c_m)=\theta(c_j)$, 
 so $c_{i_m}(|c_{i_j}|)=0$
 and $\theta(c_{i_m})=\theta(c_{i_j})$, by (c) and  (e).
Therefore, 
 $p_{i_j}\prec p_{i_m}$
and $p_{i_m}<_{\mathrm{lex}} p_{i_j}$ in $U_m$.
By the induction hypothesis, since 
 $p_j$ is $\prec$-maximal in $T_{m-1}$, $p_{i_j}$ is $\prec$-maximal in $U_{m-1}$.
Thus, 
$p_{i_m}$ is the $\prec$-immediate successor of $p_{i_j}$ in $U_m$
and $p_{i_m}<_{\mathrm{lex}} p_{i_j}$.
Hence, $U_m$ is isomorphic to $T_m$.

\underline{Case 2}. 
$p_m$  is a $\prec$-maximal node  in $T_m$
and   $p_j<_{\mathrm{lex}} p_m$, where 
 $p_j$ is the immediate $\prec$-predecessor of $p_m$ in $T_m$.
That is, $p_m$  starts a new ray over $T_{m-1}$.
Note that  
$c_m(|c_j|)=1$ so  also $c_{i_m}(|c_{i_j}|)=1$ by (c),
 which implies that
$p_{i_j}\prec p_{i_m}$ and   $p_{i_j}<_{\mathrm{lex}} p_{i_m}$.
By (e), $\theta(c_m)\not \in \{\theta(c_\ell):\ell<m\}$ implies  $\theta(c_{i_m})\not\in\{\theta(c_{i_\ell}):\ell<m\}$.
Thus, $p_{i_m}$ starts a new ray over $U_{m-1}$, 
and so $p_{i_m}$  is a $\prec$-maximal node in $U_m$.
To show that $U_m$ is isomorphic to $T_m$, it remains to show that $p_{i_m}$ is an immediate $\prec$-successor of $p_{i_j}$ in $U_m$.

By our enumeration of the nodes in $\Psi^*$, $p_j$ is not maximal in $T_{m-1}$; let $k<m$ be such that $p_k$ is the immediate $\prec$-successor of $p_j$ satisfying $p_k<_{\mathrm{lex}}p_j$.
Then $p_j=p_k\wedge p_m$.
By the induction hypothesis, $p_{i_k}<_{\mathrm{lex}} p_{i_j}$ and $p_{i_k}$ is the  immediate $\prec$-successor of  $p_{i_j}$ in $U_{m-1}$. 
 Hence,    $p_{i_j}=p_{i_k}\wedge p_{i_m}$ in $\Psi^*$. 
It follows that  $p_{i_j}$ is the immediate  $\prec$-predecessor of $p_{i_m}$ in $U_m$,
since $p_{i_m}$  starts a new ray in $U_m$.
Hence, $U_m$ is isomorphic to $T_m$.
(If we had used an enumeration of $\Psi^*$ in which $p_j$ could be $\prec$-maximal in $T_{m-1}$, 
we would just add another case in this proof.)

\underline{Case 3}. 
$p_m$ is the $\prec$-minimum  in $T_{m}$.
Let $j<m$ be such that
$p_m$ is the
 immediate $\prec$-predecessor of $p_j$ in $T_m$.
Then 
 $p_j$ is  the $\prec$-minimum  in $T_{m-1}$, so
$p_{i_j}$ is the  $\prec$-minimum in $U_{m-1}$, by the induction hypothesis.
By (e),
$\theta(c_m)=\theta(c_j)$ and $c_m(|c_j|)=2$ imply that 
$\theta(c_{i_m})=\theta(c_{i_j})$ and 
$c_{i_m}(|c_{i_j}|)=2$.
It follows that $p_{i_m}$ is the $\prec$-predecessor of $p_{i_j}$ in $U_m$.
Hence, $p_{i_m}$ is the $\prec$-minimum  in $U_{m}$; so $U_m$ is isomorphic to $T_m$.

\underline{Case 4}. 
$p_m$ is between two nodes on the same ray in $T_{m-1}$.
This means there are 
 $j,k<m$ such that 
$p_j\prec p_m\prec p_k$, where $p_k$ is the 
immediate $\prec$-successor of $p_k$ in $T_{m-1}$,
and  $\theta(c_j)=\theta(c_m)=\theta(c_k)$.
By the induction hypothesis,
 $p_{i_k}$ is the immediate $\prec$-successor of $p_{i_j}$ in $U_{m-1}$.
Further, $c_m(|c_j|)=0$ and $c_m(|c_k|)=2$,
so  $c_{i_m}(|c_{i_j}|)=0$ and $c_{i_m}(|c_{i_k}|)=2$, by (c).
Since  (e) implies that
 $\theta(c_{i_j})=\theta(c_{i_m})=\theta(c_{i_k})$, it follows that 
 $p_{i_j}\prec p_{i_m}\prec p_{i_k}$ and these three nodes are all on the same ray in $\Psi^*$ and hence in $U_m$.
Thus, $U_m$ is isomorphic to $T_m$.

In each of the four cases, we see that $U_m$ is isomorphic to $T_m$, as structures in $\mathcal{T}^*$.
By induction on $m<\om$,  the observation holds.
\end{proof}

\begin{rem}
We point out that the proof only used (c), (e), and (f) of Definition \ref{defn.ctiso}, and the fact that $f$ is one-to-one.
\end{rem}

We now define  the collection  of subtrees of $\bS$ over which we will prove Ramsey theorems that are used for the finite big Ramsey degrees of finite chains in $\Psi^*$.
Our special space of trees helps us ensure that we preserve subcopies of $\Psi^*$ obtain Ramsey theorems while preserving subcopies of $\Psi^*$. 
Further, they set us up for future work obtaining topological Ramsey spaces of such coding trees.
We follow notation in Todorcevic's book, {\em Introduction to Ramsey spaces} \cite{TodorcevicBK10}.

\begin{defn}[The space $(\mathcal{S},\le,r)$ of coding trees]\label{defn.tRs}
Let  $\langle p_n:n<\om\rangle$ enumerate $\Psi^*$ and 
 $\bS=\bS(\Psi^*)$ as in Definition \ref{defn.codingtree1}.
This $\bS$  is the maximum element of our space $\mathcal{S}=\mathcal{S}(\bS)$.
A subtree $S\sse \bS$ is a member of $\mathcal{S}$ iff $S$ is isomorphic to $\bS$.
It follows that any two $S,T\in \mathcal{S}$ are isomorphic. 
For $S\in\mathcal{S}$, 
let $\langle c^S_n:n<\om\rangle$ enumerate the coding nodes in $S$ in order of increasing length.
Equivalently, given the coding tree isomorphism $f:\bS\ra S$  we let $c^S_n$ denote $f(c_n)$.
The levels of $S$ are determined by the lengths of the coding nodes in $S$, since $S\cong \bS$ and every level of $\bS$ has a coding node.
Define the partial order  $\le$ on $\mathcal {S}$ by
$T\le S$ iff $T\sse S$.
For  $n<\om$, define the {\em $n$-th approximation} of $S$, denoted by $r_n(S)$, to be the set of the first $n$ levels of $S$;
precisely, $r_n(S)=\{s\in S:|s|< |c^S_n| \}$.
Let 
$\mathcal{AS}_n=\{r_n(S):S\in\mathcal{S}\}$ and $\mathcal{AS}=\bigcup_{n<\om}\mathcal{AS}_n$.
For $S\in\mathcal{S}$,
$\mathcal{AS}_n|S$ denotes $\{r_n(T):T\le S\}$ and $\mathcal{AS}|S$ denotes $\{r_n(T):T\le S$ and $n<\om\}$.

For $a\in \mathcal{AS}_n$ and $S\in\mathcal{S}$, 
define
\begin{align}
[a,S]&=\{T\in\mathcal{S}:r_n(T)=a\,\wedge\, T\le S\}\cr
r_{n+1}[a,S]&=\{r_{n+1}(T): T\in [a,S]\}\cr
[n,S]&=[b,S] \mathrm{\ \  where \ }b=r_n(S)
\end{align}
Define 
the transitive relation $\le_{\mathrm{fin}}$ on $\mathcal{AS}$ by 
$a\le_{\mathrm{fin}} b$ iff 
$a\sse b$ and $\max(a)\sse\max(b)$.
Define $a\le_{\mathrm{fin}} T$
iff  $a\sse T$.
Let $\depth_S(a)$ denote the least $n$ such that $a\sse r_n(S)$ if $a\le_{\mathrm{fin}} S$, and $\infty$ otherwise.

For $S\in\mathcal{S}$, define
\begin{equation}
\Psi^*|_S=\{p_n\in\Psi^*: c^{\bS}_n\in S\}.
\end{equation}
Then $\Psi^*|_S$ is a subcopy of the pseudotree $\Psi^*$, and its inherited enumeration produces a generic sequence isomorphic to $\langle T_n:n<\om\rangle$.
That is,
$\Psi^*$ restricted to the first $n$ coding nodes of $S$ is isomorphic to $T_n$.
\end{defn}

\begin{rem}
Note that (d) of Definition \ref{defn.ctiso}
implies that for $S\in\mathcal{S}$,
 for any coding node $c^S_n$ in $S$, 
 the immediate successor $t$ of $c^S_n$  with $t(|c^S_n|)=1$
 will have  $\theta$ value strictly greater than $\theta(c^S_k)$, where $c^S_k$ is the least coding node in $S$ extending $c^S_n$.
 This corresponds to  a new ray 
 branching off of $\Psi^*$ from the $n$-th node $p^S_n$ in the enumeration of $\Psi^*|_S$ which possibly skips over some
  rays in $\Psi^*$.
\end{rem}

Let $\mathcal{T}$ denote the set of all subtrees $T\sse\bS$ such that
\begin{enumerate}
    \item[(a)]
    $T$ has no terminal nodes;
\item[(b)]
the coding nodes in $T$ are exactly the nodes which are triply branching in $T$;
\item[(c)]
each node in $T$ is either a coding node, or has  exactly one immediate successor extending either by $0$ or $2$;
\item[(d)] 
for each coding node $c$ in $T$, the set of nodes  in $T$ extending $c$ by only $0$'s and $2$'s  forms a perfect subtree $T'\sse T$, and the set of coding nodes in this 
subtree $T'$ are dense in $T'$;
\item[(e)] if $c'$ is a coding node in $T$, $c$ is the maximal coding node in $T$  such that $c\subset c'$, and either $c^{\frown}0 \sse c'$ or $c^{\frown}2 \sse c'$, then $\theta(c')=\theta(c)$.
\end{enumerate}
Then each $T\in \mathcal{T}$ represents a subset of $\Psi^*$ that is again a pseudotree.
Further, $\mathcal{S}\sse \mathcal{T}$, and  each $T\in\mathcal{T}$ contains an $S\in\mathcal{S}$.
We will use this fact in proofs.
However, we work in $\mathcal{S}(\bS)$ as it helps ensure that our induction proofs actually build trees that represent pseudotrees; additionally they set the stage for a forthcoming paper obtaining topological Ramsey spaces that directly recover the big Ramsey degrees.
The reader  may fix an $\bS$  of a particular form that they like and work in that $\mathcal{S}(\bS)$.


\section{Diaries and outline of the proof of the Main Theorem}\label{sec.diaries}

This section introduces diaries for finite chains in $\Psi$ and provides an overview of the main points of the rest of the paper.  We hope this will convey intuition and  provide guideposts for the reader as they proceed through the technicalities in Sections \ref{sec.HLVariant} and \ref{sec.mainresults}.

\subsection{Diaries}\label{subsec.Diaries}

In this subsection we define diaries, diary embeddings, and diary-shaped subsets of $\bS$, 
and point out some basic observations about them.
The aim of diaries is  to record, as simply as possible, the minimal collection of unavoidable patterns arising from finite chains in $\Psi$.
The next notion of {\em almost antichain} is a step towards that.

\begin{defn}\label{defn.aac}
We call a  meet-closed subset $A\sse 3^{<\om}$  with coding nodes an {\em almost antichain} iff the following hold:
Each pair of coding nodes $c\ne d$ in $A$ has different lengths and  is either 
incomparable 
or else  $c^{\frown}1\sse d$ or $d^{\frown}1\sse c$.
If $s$ is the meet of  two incomparable nodes in $A$, then  each node in $A$ extending $s$ 
must extend either
 $s^{\frown}0$ or $s^{\frown}2$. 
\end{defn}

 The splitting nodes in an almost antichain $A$ are not considered coding nodes in $A$,
and any coding node in $A$ is either terminal or extends by $1$.
 Lemma \ref{lem.existsaac} will show that for  each $S\in\mathcal{S}$ there is 
an almost antichain $A\sse S$  such that the set 
\begin{equation}\label{eq.P}
\Psi^*|_A:=\Psi^*\re \{p_i:c_i\mathrm{\ is\ a\ coding\ node\ in\  } A\}
\end{equation}
 is isomorphic  to  $\Psi^*$.
The next observation then follows immediately from  Observation \ref{obs.easy}.

\begin{observation}\label{obs.chaininA}
Let $A\sse \bS$ be an almost antichain such that $\Psi^*|_A$ is isomorphic to $\Psi^*$.
Suppose $p_i,p_j\in \Psi^*|_A$, $i\ne j$, and $p_i$ and $p_j$ are $\prec$-comparable. There are three possibilities:
\begin{enumerate}
\item[$\bullet$]
$c_j\sse c_i$.  In this case,
 ${c_j}^{\frown}1\sse c_i$,  $p_j\prec p_i$, $p_j<_{\mathrm{lex}} p_i$, and  $p_i$ and $p_j$  are on different rays in $\Psi^*|_A$.
\end{enumerate}
Otherwise, $c_i$ and $c_j$ are incomparable and (w.l.o.g.) suppose $c_i<_{\mathrm{lex}}c_j$.
Let $k$ be such that $c_k=c_i\cap c_j$ in $\bS$.
Then $c_k$ is a splitting (but not coding) node in $A$ and 
$p_k\not\in \Psi^*|_A$.
There are now two possibilities:
\begin{enumerate}
\item[$\bullet$]
 ${c_k}^{\frown}1\sse c_i$ and $c_j= {c_k}^{\frown}2^{\frown}s$ where $s\in \{0,2\}^{<\om}$.
In this case  $p_j\prec p_i$,  $p_i<_{\mathrm{lex}}p_j$, and   they are on different rays.

\item[$\bullet$]
 ${c_k}^{\frown}0\sse c_i$ and $c_j={c_k}^{\frown}2^{\frown}s$ where $s\in \{0,2\}^{<\om}$.
In this case,
$p_j\prec p_i$ and $p_i<_{\mathrm{lex}}p_j$,
and
 $p_i$ and $p_j$ are in the same ray iff
$c_i=  {c_k}^{\frown}0^{\frown}t$ where $t\in \{0,2\}^{<\om}$.
\end{enumerate}
\end{observation}

At this point,  we are ready to  define diaries for chains.  These extract the essential information about how a chain in $\Psi^*$ can lie in $\bS$. 
We only define diaries for chains as the work in \cite{CEW_EUROCOMB25} showed that antichains in $\Psi$ of size two or more have infinite big Ramsey degrees.

\begin{defn}[Diaries]\label{def.diary}
A finite subtree $\Delta\sse 3^{<\om}$ is a {\em diary} for a finite chain in $\Psi$
iff the following hold:
\begin{enumerate}
\item
Each level of $\Delta$ has
exactly one node of only one of the following types:
\begin{enumerate}
\item[(a)]
a splitting node $s$ with exactly two immediate successors, $s^{\frown}0$ and $s^{\frown}2$;
\item[(b)]
a coding node which is terminal in $\Delta$;
\item[(c)]
a coding node  $c$ which is not terminal in $\Delta$ and has exactly one immediate successor, $c^{\frown}1$;
\item[(d)]
a non-coding node $t$ which has exactly one immediate successor, $t^{\frown}1$.
\end{enumerate}
We call these nodes  of types (a)--(d) the {\em critical nodes}  of $\Delta$ and denote them by $\langle d_j:j<h\rangle$, where
$|d_j|=j$ and  $h=\height(\Delta):=1+\max\{|t|:t\in\Delta\}$ is the {\em height} of the diary.
\item
 $\Delta\sse 3^{<h}$ and is closed under initial segments.
\item
Critical nodes of type (c) and (d) can only occur on the leftmost branch of $\Delta$.
\item
A critical node of type (d)  never occurs as the the minimal node in $\Delta$, nor as the  immediate successor of a critical node of type (c) or (d). 
\item
If $t$ is a node in $\Delta$ which is not a critical node, then  exactly one of $t^{\frown}0$ or  $t^{\frown}2$  is its immediate successor in $\Delta$.
\end{enumerate}
In particular, a diary is an almost antichain.
\end{defn}

\begin{rem}\label{rem.Devlin}
Diaries for chains in $\Psi$ are an augmentation of the diaries for the rationals, also known as Devlin types.
Devlin types are equivlent to the subcollection of our 
diaries $\Delta$ contained in $\{0,2\}^{<\om}$.
In this case,
 all levels have critical nodes of types (a) or (b),
 (3) and (4) are irrelevant, and (5) is restricted to extensions by $0$ (never $2$).
If 
$C$ is a chain contained  in  the ray  $\NR_{\theta(p_j)}$, where $\theta(p_j)\not \in\{\theta(p_i):i<j\}$, then
$C$ is encoded by  a set of coding nodes   in $\bS$  each of which is of the form  ${c_j}^{\frown}t$  where  $t\in \{0,2\}^{<\om}$.
Hence, Devlin's result \cite{DevlinThesis} for the rationals applies (see  also \cite{CDP1} and \cite{CDP2} for the coding tree version of his result).
Any collection of  coding nodes in the same ray 
in an almost antichain 
 forms a Devlin type.
\end{rem}

\begin{observation}\label{obs.finitelymanydiaries}
For each $n\ge 1$, there are finitely many diaries for chains in $\Psi$ of length $n$.
\end{observation}

This can be seen directly from Definition \ref{def.diary}.
Fix $n\ge 1$.  A diary $\Delta$ for a chain of length $n$ has  $n$ coding nodes and at most $n-1$ many splitting nodes.
Thus, the number of levels  in $\Delta$ with nodes of types (a)--(c) is at most $2n-1$.
By (4), critical 
nodes of type (d)  can never occur in immediate succession, so they can only occur strictly between two levels in $\Delta$ with 
ciritcal nodes of types (a)--(c).
Thus, the number of levels in $\Delta$ is bounded above by  $(2n-1)+(2n-2)=4n-3$.
So the number of diaries for chains of length $n$  is 
bounded above by
the number of subtrees of  $2^{<4n-3}$, 
though the actual number is much smaller for $n\ge 2$.

We next point out how
Observation \ref{obs.chaininA}
leads to the necessity of  (3) in  the definition of diary.

\begin{prop}\label{prop.no01splits}
Let $A$ be an almost antichain representing a copy of $\Psi^*$.
If $\mathcal{C}\sse \Psi^*|_A$ is a finite  chain, then any occurrences of $1$'s in the downward  closure $\widehat{C}\sse\bS$ of the set 
$C:=\{c_i:p_i\in \mathcal{C}\}\sse  A$
can only occur in the leftmost branch of $\widehat{C}$.
\end{prop}

\begin{proof}
Suppose that 
 $c_i,c_j\in C$ are incomparable.
Let $k$ be such that $c_k=c_i\cap c_j$, and 
without loss of generality, suppose $c_i$ is left of $c_j$.
It suffices to show 
\begin{equation}\label{eq.neone}
 c_j(\ell)\ne 1,\ \forall\ |c_k|\le \ell <|c_j|
\end{equation}
for then only those $c_i$ in the leftmost branch of $\widehat{C}$ can possibly have $c_i(\ell)=1$ for some $\ell<|c_i|$.
Since $p_i$ and $p_j$ are $\prec$-comparable in $\Psi^*$,
 Observation \ref{obs.chaininA} implies that 
either
(a)
 ${c_k}^{\frown}1\sse c_i$ and $c_j= {c_k}^{\frown}2^{\frown}s$ where $s\in \{0,2\}^{<\om}$;
or
(b)
 ${c_k}^{\frown}0\sse c_i$ and $c_j={c_k}^{\frown}2^{\frown}s$ where $s\in \{0,2\}^{<\om}$.
Both cases imply equation (\ref{eq.neone}).
\end{proof}

Next we define the notion of an embedding of a diary into $\bS$ and related notions.

\begin{defn}\label{defn.DE}
Let $\Delta$ be a diary for a finite chain in $\Psi$, and list the critical nodes of $\Delta$ in order of increasing length as $\langle d_j:j<h\rangle$.
Let $u,v$ represent members of $\Delta$.
A {\em diary embedding} is an injection  $e:\Delta\ra\bS$ satisfying the following, where
$e[\Delta]$  denotes the $e$-image of $\Delta$:
\begin{enumerate}

\item
$e$ preserves extensions: $u\sse v$ iff $e(u)\sse e(v)$.

\item
$e$  preserves relative lengths:  $|u|\le |v|$ iff  $|e(u)|\le |e(v)|$.

\item
If $v$ is an immediate successor of $u$,
then $v=u^{\frown}1$ iff $\theta(e(v))\ne \theta(e(u))$.

\item
If $v$ is an immediate successor of $u$,
then 
$e(v)(|e(u)|)=v(|u|)$.

\item
$e$ preserves the type of the critical nodes:
\begin{enumerate}
\item[(a)]
$d_j$  is a  node of type (a) iff $e(d_j)$ is a splitting node in $e[\Delta]$.
\item[(b)]
$d_j$ is a node of type (b)  iff $e(d_j)$ is a coding node which is terminal in $e[\Delta]$. 
\item[(c)]
$d_j$ is a   node of type (c) iff $e(d_j)$ is a non-terminal coding node in $e[\Delta]$ with one immediate successor in $e[\Delta]$.
\item[(d)]
 $d_j$ is a node of type (d) iff 
$e(d_j)$ is  a non-splitting node that is not considered as a coding node in $e[\Delta]$ and has one immediate successor in $e[\Delta]$.
\end{enumerate}
In cases (c) and (d), (3) and (4)  imply that 
$e(v)(|e(d_j)|)=1$, where $v$ is the immediate successor of $d_j$ in $\Delta$.
For each $i\in\{a,b,c,d\}$, we call  $e(d_j)$ a critical node of type (i)  in $e[\Delta]$ iff $d_j$ is a critical node of type (i) in $\Delta$.
\end{enumerate}
\end{defn}

It follows from (1)--(5) that $e$ also preserves meets:  $e(u\cap v)=e(u)\cap e(v)$.
Let $\Emb(\Delta,\bS)$ denote the set of all diary embeddings of $\Delta$ into $\bS$.
A set of  nodes  $D\sse\bS$  is  called {\em diary-shaped} iff $D=e[\Delta]$ for some  diary embedding $e$ of some  diary $\Delta$.
Given such a $D$, 
let  
$I(D)$ denote the set of those $ i\in\om$ such that  $c_i=e(d_j)$, where $d_j$ is of type (b) or (c) and
$j<h(\Delta)$.
The
{\em substructure of $\Psi^*$ encoded  by $D$ }  is defined to be
\begin{equation}\label{eqn.CNeDelta}
\Psi^*|_D= \{p_i: i\in I(D)\}.
\end{equation}
We say that two diary-shaped sets $D,D'\sse \bS$ 
are {\em diary-isomorphic} iff 
there is a diary $\Delta$ and diary embeddings $e:\Delta\ra D$ and $e':\Delta\ra D'$.

\subsection{Outline of proof of Main Theorem}

The  proof of the Main Theorem can be broken into three main ideas:
\begin{enumerate}
\item[I.]
Every diary for a finite chain in $\Psi$ has the Ramsey property. 
\item[II.]
There is a subcopy of  $\Psi^*$ in which every finite chain is diary-shaped.
\item[III.]
For each $n$, there are finitely many diaries 
for chains of size $n$ in $\Psi$.
\end{enumerate}

The proof of  I.\  is the subject of 
 Section \ref{sec.HLVariant}, the main result of which 
 appears as the  last sentence in Theorem \ref{thm.FBRD}.
This is a  Milliken-style  Theorem for diaries.
   II.\ will be proved in  Section  \ref{sec.mainresults}, and 
III.\ has already been shown in Observation \ref{obs.finitelymanydiaries}.

In what follows, $\Delta|_m$ denotes the restriction of $\Delta$ to its first $m$ levels.

\begin{thmFBRD}
Let $\Delta$ be a diary for a finite chain in $\Psi$, and fix  $0\le m<\height(\Delta)$.
For any coloring $\chi$ of $\Emb(\Delta|_{m+1},\bS)$  into finitely many colors, there is a $\phi\in\Emb(\bS,\bS)$ 
with $\{\phi\circ g:g\in   \Emb(\Delta|_{m+1},\bS)  \}$ monochromatic for $\chi$.
In particular,  if $\chi$ is a finite coloring of $\Emb(\Delta,\bS)$, then there is a $\phi\in\Emb(\bS,\bS)$ with $\chi$ constant on $\phi\circ \Emb(\Delta,\bS)$.
\end{thmFBRD}

Theorem \ref{thm.FBRD} will be proved  
by induction on $m$  applying the following theorem.
Here,  $\tilde{f}(m-1)$ denotes  the length of the longest node(s) in $f[\Delta|_m]$.

\begin{lemEnd-homog}
Suppose $\Delta$ is a diary for a finite chain in $\Psi$ and  $m<\height(\Delta)$.  Then for any coloring $\chi$ of  $\Emb(\Delta|_{m+1},\bS)$, there is a 
 $\phi\in\Emb(\bS,\bS)$ such that 
for each $f\in\Emb(\Delta|_m,\bS)$, $\chi$ is monochromatic on 
$\{\phi\circ g:g\in G_f\}$, where 
$G_f:=\{g\in\Emb(\Delta|_{m+1},\bS)  :g|_m=f\}$.
\end{lemEnd-homog}

 Section \ref{sec.HLVariant} begins by proving some  necessary machinery in 
Subsection \ref{subsec.AL} 
 in the form of   Amalgamation Lemmas
\ref{lem.generalamalg} and 
\ref{lem.A.3}, the second of which shows that  Axiom A.3 for topological Ramsey spaces  (see \cite{TodorcevicBK10}) holds in $\mathcal{S}$.
The main result in Subsection \ref{subsec.HLVariant}, 
is
Theorem 
\ref{thm.HLVariant},  
 a new 
 variant of the Halpern--\Lauchli\ Theorem.
This 
variant  cannot be  proved  by only applying the
 Halpern--\Lauchli\ Theorem, because 
diaries can  have critical nodes of types (c) and (d), and these  are allowed to range over $\om$ many different rays.
Instead, we first prove 
 Lemma
\ref{lem.Stepone} by induction applying  the Halpern--\Lauchli\ Theorem  $\om$ many times. 
To homogenize over the different rays, we use forcing to  obtain  Theorem 
\ref{thm.HLVariant}.
In Subsection \ref{subsec.DiaryMilliken}, we 
prove 
Lemma  \ref{lem.End-homog}
 by inductively  applying Corollary \ref{cor.End-homog}, an immediate consequence of 
Theorem 
\ref{thm.HLVariant}.  Lemma  \ref{lem.End-homog}
 then forms the pigeonhole from which we prove the 
 Ramsey Theorem for diaries, Theorem \ref{thm.FBRD}.

 Section \ref{sec.mainresults} is concerned with proving the Main Theorem \ref{thmmain}.
In Subsection \ref{section.aa} we first prove the previously mentioned 
Lemma   \ref{lem.existsaac}, constructing within any $S\in\mathcal{S}$ an almost antichain $A$ representing copies of $\Psi^*$.
Then Theorem
\ref{thm.Asubsetsdiaries} shows
that  to  each finite chain $C$ in $\Psi^*|_A$ 
there corresponds a unique
 diary-shaped subset of $A$, and hence corresponds to a unique diary and diary embedding. 
The  proof of  Main Theorem then follows quickly.
Lastly, in Subsection \ref{subsec.ex7}
we display the seven diaries for chains of length two in $\Psi$.


\section{A Halpern--\Lauchli\ Variant and a Milliken Theorem for diaries }\label{sec.HLVariant}

Subsection \ref{subsec.AL} contains amalgamation lemmas that will be used throughout the rest of the paper.
In Subsection \ref{subsec.HLVariant}, we  prove a variant of the Halpern--\Lauchli\  Theorem (Theorem \ref{thm.HLVariant}); this is a Ramsey theorem 
 for finite colorings of certain level sets in coding trees.  This will be applied  in  Subsection \ref{subsec.DiaryMilliken}  to obtain Lemma  \ref{lem.End-homog}, which forms  the pigeonhole principle used in the inductive proof of Theorem \ref{thm.FBRD}, a Milliken variant for finite colorings of embeddings of a given diary into a coding tree.

We chose to present Subsections 
\ref{subsec.AL} and \ref{subsec.HLVariant}  in more generality than is necessary for the results in this paper for two reasons. Firstly, the general theorems in these sections will be applied in forthcoming work on topological Ramsey spaces consisting of certain copies of $\Psi^*$ (not defined in this paper).
Secondly, the proofs are exactly the same as for the more restricted  case of  finite diaries.

Some readers might have preferred us to use trees that are fully $3$-branching at each level.  
We have checked that the results in this section can be proved using such trees and applying a certain  one-to-one correspondence between such trees and certain trees in the class $\mathcal{T}$.
However, for the sake of stating theorems for coding trees induced by an $\om$-ordering of the nodes in $\Psi^*$ and for  applications to  future work, we prefer to work with the coding trees in $\mathcal{S}$ as defined in  Section \ref{sec.basics}.

\subsection{Amalgamation Lemmas}\label{subsec.AL}

The amalgamation lemmas in this subsection are  necessary for the work in following subsections.

\begin{defn}\label{defn.perfect}
Let $U$ be a meet-closed subset of $\bS$, (equivalently, a  meet-closed subset of any $S\in\mathcal{S}$).
A node $s\in U$  is called a {\em splitting node} in $U$
iff $s$ is the meet of two incomparable nodes in $U$.
When we that say $U$ is a {\em subtree},  it is implied that every splitting node in $U$ is also considered as a  coding node in $U$.
We call $U$ a {\em \perfect\ tree} 
  iff 
$U$ is a subtree of $\mathbb{S}$   such that 
 $U\sse \{s\in \bS:\theta(s)=k\}$ for some $k\in\om$, and 
 each node in $U$ has two incomparable extensions in $U$.
\end{defn}

Note that $U$ is a ray tree iff $\Psi^*|_U$, the substructure of $\Psi^*$ induced by the coding nodes  in $U$, is a dense linear order without endpoints that is contained in a ray $\NR_k$, for some $k$.

We will use the following notation throughout.
Given a finite set of nodes $a\sse\bS$, 
let $\max a$ denote the set
of terminal nodes in $a$ and 
 $(\max a)^+$ denote the set 
of immediate successors 
in $\bS$  of the nodes in $\max a$.
For $S\in\mathcal{S}$, let $\widehat{S}$ denote the set of all nodes in $\bS$  which are  initial segments of members of $S$.
Given $x\in \widehat{S}$, let $S_x$ denote the subtree of $S$ consisting of all nodes in $S$ extending $x$.

\begin{defn}\label{defn.leveltrees}
Suppose that $S\in \mathcal{S}$, $X\sse S$ is a set of nodes, each of the same length, and $L\sse\om$ is infinite.
We  say that a collection of subtrees $U_x\sse S_x$, $x\in X$, {\em has levels in $L$} iff
for each coding node  $c\in \bigcup_{x\in X} U_x$,  $|c|\in L$ and 
$s\re |c|\in U_x$ for each $x\in X$ and each $s\in U_x$ with $|s|\ge |c|$.
\end{defn}

\begin{lem}[Amalgamation]\label{lem.generalamalg}
Let $S\in \mathcal{S}$ and $d\in\om$ be given. Let $X$  be a non-empty set  such that   $X\sse (\max r_d(S))^+$ if $d\ge 1$,  or  $X=\{c^S_0\}$ if $d=0$.
For  $x\in X$, 
let  $U_x$ be a  subtree of $S_x:=\{s\in S:x\sse s\}$ where,
letting $c$ denote the maximum coding node in $S$ such that $c\sse x$,
one of the following hold:
\begin{enumerate}
    \item[(0)]
     $x(|c|)\in\{0,2\}$ or $d=0$,  and 
$U_x\sse \{s\in S: s\contains x\, \wedge\, \theta(s)=\theta(c)\}$ is a \perfect\ tree;
    \item[(1)] 
$x(|c|)=1$ or $d=0$, and 
$U_x\sse \{s\in S: s\contains x\, \wedge\, \theta(s)=k\}$ is a \perfect\ tree, where $k=\theta(c')$ for some $c'$ in $U$ extending $x$;
\item[(2)]  
      $x(|c|)=1$ or $d=0$, 
 all splitting nodes in $U_x$ have branching degree $3$, and $U_x$ contains a member of $\mathcal{S}$.
\end{enumerate}
For $j<3$, let $X_j$ denote the set of those $x\in X$ for which we are in Case (j).
 Then there is a $T\in [d,S]$ such that
 \begin{enumerate}
     \item[(0)]
For each $x\in X_0$, 
$\{t\in T:t\contains x\,\wedge\, \theta(t)=\theta(x)\}$ is a subtree of $\widehat{U}_x$;
\item[(1)]
For each $x\in X_1$,  letting $c_x$ denote the least coding node in $T$ extending $x$,
we have 
 $\{t\in T:t\contains x\wedge \theta(t)=\theta(c_x)\}
 =\{t\in T:t\contains c_x\wedge \theta(t)=\theta(c_x)\}$ is a subtree of 
$\widehat{U}_x$.
\item[(2)]
For each $x\in X_2$, 
 $T_x\sse \widehat{U}_x$ (and  $T_x$ contains a member of $\mathcal{S}$).
 \end{enumerate}
Further, if the $U_x$, $x\in X$, have levels in 
 an infinite set $L\sse \om$,
 then  $T$ can be chosen so that for each $x\in X$, $T_x\sse U_x$.
\end{lem}

\begin{proof}
The idea of the proof is simple:  Build a tree $T$ extending  $r_d(S)$ inside $S$ by
choosing any coding nodes that can be in a $U_x$ to be in a $U_x$.  More generally,
whenever you can extend and stay inside one of the $\widehat{U}_x$'s, do so; otherwise, just build inside $S$ as needed to obtain a member of $[d,S]$.

For $x\in X_0$, let $N_x$ denote the set of all $n>d$ for which  the maximum coding node in $r_n(S)$ extends $x$ and has $\theta$-value equal to $\theta(x)$.
For $x\in X_1$, let $c^S_x$ denote the least coding node in $S$ extending $x$, and 
let $N_x$ denote the set of all $n>d$ for which the maximum coding node in $r_n(S)$ 
extends $c^S_x$ and has $\theta$-value equal to $\theta(c^S_x)$.
For $x\in X_2$, let $N_x$ denote the set of all $n>d$ for which $c^S_n$, the maximum coding node in $r_{n+1}(S)$,
extends $x$.
Let $n_0<n_1<\dots$ enumerate the set $N:=\bigcup_{x\in X}N_x$.

For $i=0$:
Let $y$ be the member of $X$ such that the maximum coding node of any member of $r_{n_0}[d,S]$  extends $y$.
Take $b_0\in r_{n_0+1}[d,S]$ so 
that  $c^{b_0}_{n_{0}}\in U_{y}$, and 
for every other  $x\in X\setminus\{y\}$,
the (unique) extension  of $x$  in $\max(b_0)$  is a member of  $\widehat{U}_x$.

Now suppose that  $i>0$ and  $b_{i-1}\in r_{n_{i-1}+1}[d,S]$ satisfies the conclusion of the lemma.
Let $j<3$  and $y\in X$
be such that $n_i\in N_{y}$, where $y\in X_j$.
Take $b_i\in  r_{n_i+1}[b_{i-1},S]$ so that
 the following hold:
The coding node  $c^{b_i}_{n_i}$ in $\max(b_i)$ is in $U_j$.
Furthermore, 
for  each  $t\in \max(b_i)\setminus \{ c^{b_i}_{n_i} \}$  that extends some $x\in X$, 
 we require:
If  $t$  extends some $x\in X_0$ and 
$\theta(t)=\theta(x)$,
then $t\in \widehat{U}_x$.
If  $t$ extends some $x\in X_2$, then   $t\in \widehat{U}_x$.
(This is possible since,
by the induction hypothesis, each node in $b_{i-1}$ extending $x$ is a member of $\widehat{U}_x$.)
Lastly, suppose  that  $t$  extends some $x\in X_1$.
If $n_i<\min(N_x)$, then make sure that the (unique) node in $\max(b_i)$ extending $x$ is a member of 
 $\widehat{U}_x$.
If $n_i=\min(N_x)$, then 
we have already ensured that 
 the coding node $c^{b_i}_{n_i}$ 
is in 
 $U_x$. 
For this special case, let $c_x^*$ denote this coding node, noting that $\theta(c_x^*)=k$, 
 and use in the next case.
If $n_i>\min(N_x)$ and the 
initial segment $s$ of $t$ in $\max(b_{i-1})$ satisfies $\theta(s)=k$, 
then by the induction hypothesis, $s$ is in $\widehat{U}_x$ and we can take $t$ to be in $\widehat{U}_x$.

Let $T=\bigcup_{i<\om} b_i$.  Then $T\in [d,S]$
and satisfies the conclusion of the lemma.
Recall  that the coding nodes  in the $\max(b_i)$ were chosen to be coding nodes in the trees $U_x$, $x\in X$.
In the case that there is an infinite $L$ where 
 for all $x\in X$, $U_x\cap 3^{\ell}=\{s\in U_x:|x|=\ell\}$ for each $\ell\in L$,
then  choose the coding nodes in the $\max(b_i)$ to be nodes in the $U_x$. Then the rest of the
 nodes $t$ chosen in the construction of the $b_i$ may be chosen to be nodes  in the trees $U_x$.
\end{proof}

The next lemma shows that Axiom \bf A.3 \rm from Todorcevic's four axioms guaranteeing a topological Ramsey space. (See Chapter 5 of  \cite{TodorcevicBK10}.)

\begin{lem}[{\bf A.3}]\label{lem.A.3}
$\mathcal{S}$ satisfies Axiom
{\bf  A.3}.
\end{lem}

\begin{proof}
Axiom {\bf A.3} states that for $a\in\mathcal{AS}$ and $S\in\mathcal{S}$,
\begin{enumerate}
    \item 
If $\depth_S(a)<\infty$ then $[a,S']\ne\emptyset$ for each $S'\in [\depth_S(a),S]$;    \item 
$S'\le S$ and $[a,S']\ne\emptyset$ imply there is a $T\in[\depth_S(a),S]$ such that 
$\emptyset\ne [a,T]\sse[a,S']$.
\end{enumerate}

Towards proving (1),  let $d=\depth_S(a)$ and let $S'\in [d,S]$ be given. 
We will show that $[a,S']\ne\emptyset$ by showing there is a $T\le S'$ which extends $a$.
Let $X=(\max a)^+$, and for each $x\in X$, let $U_x$ denote the set of nodes in $S'$ extending $x$.
By  Lemma \ref{lem.generalamalg},
there is a  $T$ end-extending $a$ 
with  $T\setminus a\sse\bigcup\{U_x:x\in X\}$.
Hence, $T$  is in $[a,S']$.

For (2), let $d=\depth_S(a)$ and 
suppose $S'\le S$  satisfies $[a,S']\ne\emptyset$.
Again, let $X=(\max a)^+$ and for each $x\in X$, let $U_x$ denote the set of nodes in $S'$ extending $x$.
By Lemma \ref{lem.generalamalg}, there is a  $T\in [d,S]$ such that for each $x\in X$, all nodes extending $x$ in $T$ are members of $U_x$, and hence in $S'$.
Thus,  $T\in [a,T]\sse [a,S']$.
\end{proof}

\subsection{A Halpern--\Lauchli\ Variant}\label{subsec.HLVariant}

First we recall the following  variation for perfect trees  of the Halpern--\Lauchli\ Theorem \cite{Halpern/Lauchli66}.
We state it in terms of ray trees, as that will be our  application.

\begin{thm}\label{thm.HLSkew}
Let $1\le n<\om$ be given and $i_*<n$ be fixed.
For each $i<n$,
let  $U_i$ be a 
 \perfect\ tree.
Let $L=\{|t|\in U_{i_*}:t$ is a splitting node in $U_{i_*}\}$.
Define 
$$
\bigotimes_{i<n, i_*} U_i=
\bigcup_{\ell\in L}\prod_{i<n} U_i\re \ell
$$
Let $c$ be any coloring of $\bigotimes_{i<n, i_*} U_i$ into finitely many colors.
Then there are \perfect\  subtrees $V_i\sse U_i$ for each $i<n$ such that $c$ is constant on $\bigotimes_{i<n, i_*} V_i:= \bigcup_{\ell\in M}\prod_{i<n} V_i\re \ell
$, where $M$ is the set of lengths of splitting nodes in $V_{i_*}$.
\end{thm}

This theorem   has been  applied in works of 
Huber--Geschke--Kojman \cite{HGK2019} and Y.Y. Zheng \cite{Zheng18} and in Dobrinen--Wang \cite{Dobrinen/Wang19} and follows easily from the strong tree version of the Halpern--\Lauchli\ Theorem (see \cite{TodorcevicBK10}).

Recall  from Subsection \ref{subsec.Diaries}
that in  diaries for chains in $\Psi$,
 only one  node per level  is allowed to change rays by the next level of the diary.
 The reason for this is that allowing two nodes in the same level  of a coding tree to change  rays simultaneously  allows for an encoding of an oscillation function, which was used in \cite{CEW_EUROCOMB25} to show that antichains of size more than one  in $\Psi$ have infinite big Ramsey degree. 
Hence, any Ramsey theorem extending the Halpern--\Lauchli\ Theorem to  our coding trees must restrict to only allowing one node to change rays at a time.
This is reflected in the   following Assumption \ref{assum.prodtrees} by the fact that the set $X_2$  (recall Case (2) in the Amalgamation Lemma \ref{lem.generalamalg}) can have at most one element.  
Recall also that  diaries for chains in $\Psi$ can only have ray changes occuring in the leftmost branch. 
So for the proof of the Main Theorem we would only need to consider that case below.  However, as the full versions  of 
Lemma \ref{lem.Stepone} and Theorem \ref{thm.HLVariant}
 will be used in a following paper, we prove these more general statements, which are no harder to prove than their restriction to diaries.

The following assumption and notation  will be used in our Ramsey theorems in this subsection.

\begin{assumption}\label{assum.prodtrees}
Suppose we are given $S\in \mathcal{S}$ and $d\in\om$. 
If $d=0$,  let  $X =\{c^S_0\}$.
If $d\ge 1$, let $X$ be a subset of $(\max r_d(S))^+$.
The following restricted version of the  hypotheses of Lemma \ref{lem.generalamalg} is assumed: 
$X=X_0\cup X_2$, where 
$X_2$ is either empty or a singleton denoted by $\{y\}$. 
For each $x\in X_0$,
$U_x\sse \{s\in S: s\contains x\, \wedge\, \theta(s)=\theta(x)  \}$ is a  \perfect\ tree.
In the case that $X_2\ne\emptyset$ and $d=0$, then $X_2=X$ and $y=c^S_0$.
In the case that $X_2\ne\emptyset$ and $d\ge 1$, then
 $y\contains {c_y}^{\frown}1\sse y$, where 
 $c_y$ is 
the maximum coding node in $S$ 
below $y$.
In both cases, 
$U_y\sse \{s\in S: s\contains y\}$ is in 
$\mathcal{T}$.

The subtrees $U_x$, $x\in X$, have levels in 
 an infinite set $L\sse \om$.
Fix  an $x_*\in X$.
Suppose $f$ colors the level product
$\bigotimes_{x\in X, x_*} U_x$
 into finitely many colors.
If $x_*\in X_2$,
  let $\langle n_i:i<\om\rangle$  denote the strictly increasing sequence of those $n\in\om$ for which
 $c^S_n\contains x_*$ and 
 $\theta(c^S_n)\ne \theta(c^S_m)$ for all $m<n$.
\end{assumption}

The following is the closest to homogeneity that we know how to obtain without using forcing.

\begin{lem}[End-Homogeneity]\label{lem.Stepone}
Let $S,d,X,X_0,X_2,y, x_*,U_x\ (x\in X),f$ be as in
Assumption  \ref{assum.prodtrees}.
Then there is a $T\in [d,S]$ such that 
the following hold:
For $x\in X_0$, 
$V_x:= \{t\in T: x\sse t \wedge \theta(t)=\theta(x)\}\sse U_x$.
  Further, 
\begin{enumerate}
\item[$\bullet$]
If  $X_2=\emptyset$, then $f$ is constant on $\bigotimes_{x\in X,x_*}V_x$.
\item[$\bullet$]
If  $X_2\ne\emptyset$,
 then $V_y:=\{t\in T:y\sse t\}\sse U_y$.
For each
$i<\om$  let $\ell_i=|c^T_{n_i}|$.
Then for each $\vec{w}\in \prod_{x\in X}V_x\re \ell_i$ 
 with $w_y=c^T_{n_i}$,
 we have 
$f$ 
  constantly equal to $f(\vec{w})$ on the set of all
$\vec{t}\in\bigotimes_{x\in X, x_*}V_x$ 
with
$t_x\contains w_x$ for  $x\in X$ and $\theta(t_y)=\theta(c_{n_i}^T)$.
\end{enumerate}
\end{lem}

\begin{proof}
Suppose first that 
$X_2=\emptyset$.
Then for  each  $x\in X$, $U_x$ is a subtree of $\{s\in S:s\contains x\, \wedge\, \theta(s)=\theta(x)\}$.
By  Theorem 
\ref{thm.HLSkew}, there are subtrees $V_x\sse U_x$ such that $f$ is constant on
$\bigotimes_{x\in X, x_*}V_x$.
Apply   Lemma \ref{lem.generalamalg} to obtain $T\in [d,S]$ so that for each $x\in X$,
$T\cap U_x\sse V_x$.
Then the conclusion of this theorem holds.

Now suppose that $X_2\ne\emptyset$ and $x_*=y$.
We will do an induction scheme where at each stage $i<\om$, we 
apply Theorem \ref{thm.HLSkew}  and  the Amalgamation Lemma successively, finitely many times.
This will construct a subtree of $S$ extending $r_d(S)$ which will represent a sub-pseudotree of $\Psi$ which we can  thin to obtain a coding tree in $[d,S]$.
\vskip.1in

\underline{Base Case}. 
Let $y_{0}$ denote the least coding node in $U_y$
and let $U_{y_0}=\{s\in U_y:\theta(s)=\theta(y_0)\}$.
Note that  $y_0=c^S_{n_0}$.
 Let $\ell=|y_0|$.
For $x\in X_0$, let 
$I_x= U_x\re\ell
:=
\{w \in U_x:  |w|=\ell\}$; 
let $I_y=\{y_{0}\}$.
 For  $\vec{w}=\langle w_x:x\in X\rangle$ in $\prod_{x\in X}I_x$, we may
apply Theorem \ref{thm.HLSkew}
  to obtain subtrees $U^{\vec{w}}_{w_x}\sse \{ t\in U_x: t\contains w_x\}$
  and $U^{\vec{w}}_{y_0}\sse U_{y_0}$ such that 
  $f$ is constant on $\bigotimes_{x\in X,y}U^{\vec{w}}_{w_x}$.
Then apply the Amalgamation Lemma \ref{lem.generalamalg} to obtain an $S^{\vec{w}}\in [n_0,S]$ such that for each $x\in X_0$,
all  nodes in $S^{\vec{w}}$ extending $w_x$ are in $U^{\vec{w}}_{w_x}$, and all nodes in $S^{\vec{w}}$ extending $y_0$ with $\theta$-value equal to $\theta(y_0)$ are in $U_{y_0}$.

Do this process successively in a finite induction over all
$\vec{w}$ in $\prod_{x\in X}I_x$
 to obtain a $T_0\in [n_0,S]$ 
 with the following properties:
For each $x\in X_0$ and $w\in I_x$, 
letting
$V_{w}:=\{s\in T_0: s\contains w_x$ and $ \theta(s)=\theta(x)\}$, we have that $V_w
\sse U_x$; and
letting 
 $V_{y_0}=\{s\in T_0:\theta(s)=\theta(y_0)\}$, we have that $V_{y_0}\sse U_y$.
Most importantly,
for each sequence $\vec{w}=\langle w_x:x\in X\rangle\in\prod_{x\in X}I_x$,
 $f$ is constant on  $\bigotimes_{x\in X,y} V_{w_x}$.
\vskip.1in

\underline{Induction Hypothesis}. 
Suppose $1\le i<\om$ and we have $T_j$ for $j< i$ such that   for each $j<i$
 (letting $T_{-1}=S$) the following hold:
 \begin{enumerate}
   \item  
$T_{j}\in[n_j,T_{j-1}]$.
\item 
Let $y_{j}$ denote the coding node $c^{T_{j-1}}_{n_j}$.
Note that $y_j$ extends $y$ and is least in $T_{j-1}$ with its $\theta$-value.
Let
 $U_{y_{j}}=\{s\in T_{j-1}:  s\contains y_{j},\ s\in U_y,$ and $\theta(s)=\theta(y_{j})\}$.
Let $\ell=|y_j|$.
For $x\in X_0$, 
let 
$I_x=\{w\in T_{j-1}\re \ell:w\contains x \,\wedge\, \theta(w)=\theta(x)\}$; 
note that $I_x\sse U_x$.
Let $I_y=\{y_j\}$.
For each sequence $\vec{w}=\langle w_x:x\in X\rangle\in\prod_{x\in X} I_x$, 
$f$ is constant on $\bigotimes_{x\in X,y} V_{w_x}$, where 
for $x\in X$, $V_{w_x}:=\{s\in T_j: 
s\contains w_x$ and $\theta(s)=\theta(x)\}\sse U_x$. 
 \end{enumerate}
\vskip.1in

\underline{Induction Step}.
Let $y_i$ denote the coding node   $c^{T_{i-1}}_{n_i}$; this is the least coding node  in $T_{i-1}$ extending $y$ above $|y_{i-1}|$.
Let
 $U_{y_i}=\{s\in T_{i-1}: s\contains y_i$, $ s\in U_y$, and $ \theta(s)=\theta(y_i)\}$.
 Let $\ell=|y_i|$.
For $x\in X_0$, let 
$I_x=\{w\in T_{i-1}\re\ell:w\contains x\ \wedge\  \theta(w)=\theta(x)\}$; 
let $I_y=\{y_i\}$.
Just as in the Base Case, 
do a finite induction over all $\vec{w}\in \prod_{x\in X} I_x$, at each step doing one application of  Theorem \ref{thm.HLSkew} and one application of the Amalgamation Lemma \ref{lem.generalamalg}  to obtain 
a $T_i\in [n_i, T_{i-1}]$ satisfying (2) of the Induction Hypothesis (with $j$ replaced by $i$).

After the $\om$ steps of this induction, let   $T=\bigcup_{i<\om}r_{n_i}(T_i)$.
Then $T\in [d,S]$ and $T$ satisfies the conclusion of the Lemma.
The proof for the case that $X_2\ne\emptyset$ but  $x_*\in X_0$ is similar.
\end{proof}

The next theorem 
does what the previous lemma could not:
In the case that $X_2\ne\emptyset$,  it homogenizes  the 
over the various different rays 
$\NR_{\theta(c)}$ for coding nodes $c$
extending $y$.
The proof uses forcing on the coding trees in $\mathcal{S}$.  
While the proof uses elements of prior proofs
 (for instance, \cite{DobrinenJML20} and \cite{CDP1}), the novelty here is in the restriction that certain extensions of certain nodes must preserve $\theta$-values.

\begin{thm}[Halpern-\Lauchli\ Variant]\label{thm.HLVariant}
Let $S,d,X,X_0,X_2,y, x_*,U_x\ (x\in X),f$ be as in
Assumption  \ref{assum.prodtrees}.
Then there is a $T\in [d,S]$ such that 
$f$ is constant on $\bigotimes_{x\in X,x_*}V_x$, where for $x\in X_0$, 
$$
V_x:= \{t\in T: x\sse t \wedge \theta(t)=\theta(x)\}\sse U_x,
$$
and if  $X_2\ne\emptyset$, then 
$$
V_y:=\{t\in T:y\sse t\}\sse U_y.
$$ 
\end{thm}

\begin{proof}
We may  assume that we have already applied  Lemma \ref{lem.Stepone}. 
Let $U_y'$ denote the set of those coding nodes $c\in U_y$ such that $c$ is least with its $\theta$-value in $U_y$.
Note that the downwards closure of $U_y'$ contains $U_y$.
Let $X^*=X\setminus \{x_*\}$.
For $x\in X^*$, let $U'_x=U_x$.
If $x_*\in X_0$, let  $M$ denote the set of those  $m<\om$  so that $\{c^S_m: m\in M\}$ is the set   of all  coding nodes extending $x_*$ in $U_{x_*}$ with $\theta(c^S_m)=\theta(x_*)$, and let $E$ denote $U'_{x_*}$.
If $x_*=y$, let   $M$ denote the set of those  $m<\om$  such that  $c^S_m\in U_y'$,
and let $E$ denote $U_y'$.

Let $k=|X^*|$ and 
let $\kappa$ be  large enough so that $\kappa\ra (\aleph_1)^{2k}_{\aleph_0}$.
By the \Erdos-Rado Therem, $
\kappa=\beth_{2k-1}(\aleph_0)^+$ suffices.
Let $L$ be the set of lengths of the coding nodes in $E$.
Define  $\bP$ to be the set of 
conditions $p$ satisfying the following:
There is a  $\vec{\delta}_p\in [\kappa]^{<\om}$ and an 
$\ell_p\in\om$ so that 
$p$ is a function with 
domain 
$\{x_*\}\cup (X^* \times\vec{\delta}_p)$ and  the following hold:
\begin{enumerate}
\item
For   $(x,\delta)\in X^*\times \vec{\delta}_p$,
$x\sse  p(x,\delta)\in U'_x\re\ell_p$;
\item
$p(x_*)$  is the coding node in $E$ of length $\ell_p$.
\end{enumerate}
The partial  ordering on $\bP$ is defined by 
$q\le p$ if and only if
\begin{enumerate}
\item
For $(x,\delta)\in X^*\times \vec{\delta}_p$, 
$p(x,\delta)\sse q(x,\delta)$;
\item
$p(x_*)\sse q(x_*)$.
\end{enumerate}
This partial ordering is atomless and separative.

Forcing with
$\bP$ adds 
one infinite branch of coding nodes in $E$ and 
$\kappa$ many branches through 
$U'_x $  for each 
 $x\in X^*$. 
 However, a generic filter for $\bP$  will not give us the theorem we seek. Rather we use the Truth Lemma $\om$ along with the Amalgamation Lemma $\om$ many times to  build the desired $T\in[d,S]$.

For  $x\in X^*$ and $\al<\kappa$,
let  $\dot{b}_{x,\al}=\{\langle p(x,\al) , p  \rangle:p\in\bP\}$, the {\em $\al$-th generic branch through $U'_x$}.
Let $\dot{b}_{x_*}=\{\langle p(x_*) , p\rangle:p\in \bP\}$, the {\em generic branch of coding nodes through $E$}.
Let $\dot{L}$ be a name for the lengths of nodes in $\dot{b}_{x_*}$, and let
 $\dot{\mathcal{U}}$ be a name for a non-principal ultrafilter on $\dot{L}$.
To simplify notation, we  assume the lexicographic order on the nodes in  $X^*$ and
write sets $\{\al_x:x\in X^*\}$ in
$[\kappa]^d$ as vectors $\vec{\al}=\lgl \al_x: x\in X^*\rgl$ in strictly increasing order.
For $\vec{\al}\in[\kappa]^d$,
we let
\begin{equation}
\dot{b}_{\vec{\al}}\mathrm{\  \  denote\ \ }
\lgl \dot{b}_{x_*}\rangle^{\frown}\lgl \dot{b}_{x,\al_x}:x\in X^*\rgl,
\end{equation}
 and for any $\ell\in \dot{L}$, let 
\begin{equation}
\dot{b}_{\vec\al}\re \ell
\mathrm{\ \ denote \  \ }
\{\dot{b}_{x_*}\re \ell\}\cup
\{\dot{b}_{x,\al_x}\re \ell:x\in X^*\}.
\end{equation}

For each $\vec\al\in[\kappa]^d$,
choose a condition $p_{\vec{\al}}\in\bP$ such that
\begin{enumerate}
\item
 $\vec{\al}\sse\vec{\delta}_{p_{\vec\al}}$;
\item
$p_{\vec{\al}}\forces$ ``There is an $\varepsilon\in 2$  such that
$f(\dot{b}_{\vec{\al}}\re \ell)=\varepsilon$
for $\dot{\mathcal{U}}$ many $\ell$";
\item
$p_{\vec{\al}}$ decides a value for $\varepsilon$, label it  $\varepsilon_{\vec{\al}}$; and
\item
$f(\{p_{\vec{\al}}(x_*)\}\cup \{p_{\vec\al}(x,\al_x):x\in X^*\})=\varepsilon_{\vec{\al}}$.
\end{enumerate}

Such conditions $p_{\vec\al}$ may be obtained as follows.
Given $\vec\al\in[\kappa]^d$,
take $p_{\vec\al}^1$ to be any condition such that $\vec\al\sse\vec{\delta}_{p_{\vec\al}^1}$.
Since $\bP$ forces $\dot{\mathcal{U}}$ to be a nonprincipal ultrafilter on $\dot{L}$, there is a condition
 $p_{\vec\al}^2\le p_{\vec\al}^1$ such that
$p_{\vec\al}^2$ forces that $f(\dot{b}_{\vec\al}\re \ell)$ is the same color for $\dot{\mathcal{U}}$ many $\ell$.
Furthermore,  there must be a stronger
condition deciding which of the colors
$f(\dot{b}_{\vec\al}\re \ell)$ takes on $\dot{\mathcal{U}}$ many levels $\ell$.
Let $p_{\vec\al}^3\le p_{\vec\al}^2$  be a condition which decides this  color, and let $\varepsilon_{\vec\al}$ denote  that color.
Finally, since $p_{\vec\al}^3$ forces that for  $\dot{\mathcal{U}}$ many $\ell$ the color
 $f(\dot{b}_{\vec\al}\re \ell)$
will equal $\varepsilon_{\vec{\al}}$,
there is some $p_{\vec\al}^4\le p_{\vec\al}^3$ which decides some level $\ell$ with length greater than the length of the nodes in $\ran(p_{\vec\al}^3)$ so that
$f(\dot{b}_{\vec\al}\re \ell)=\varepsilon_{\vec{\al}}$.
If $\ell_{p_{\vec\al}^4}<\ell$,
let $p_{\vec\al}$ be any member of $\bP$ such that
$p_{\vec\al}\le p_{\vec\al}^4$ and $\ell_{p_{\vec\al}}=\ell$.
If $\ell_{p_{\vec\al}^4}\ge \ell$,
then
let  $p_{\vec\al}\le p_{\vec\al}^3$  be a member of $\bP$ with nodes of length $\ell$ and such that 
$p_{\vec\al}(x_*)$  is the coding node in $\ran(p_{\vec\al}^4)\re\ell$ and for each $x\in X^*$, 
$p_{\vec\al}(x,\al_x)=p_{\vec\al}^4(x,\al_x)\re\ell$. 
Then $p_{\vec\al}$ forces that $\dot{b}_{\vec\al}\re \ell=\{p_{\vec\al}(x_*)\}\cup\{p_{\vec\al}(x,\al_x):x\in X^*\}$, and hence
 $p_{\vec\al}$  forces that
$f(\dot{b}_{\vec\al}\re \ell)=\varepsilon_{\vec{\al}}$.

Recall that we chose  $\kappa$ large enough so that   $\kappa\ra(\aleph_1)^{2d}_{\aleph_0}$ holds.
Now we prepare  for an application of the \Erdos-Rado Theorem.
Given two sets of ordinals $J,K$ we write $J<K$ if and only if every member of $J$ is less than every member of $K$.
Let $D_e=\{0,2,\dots,2d-2\}$ and  $D_o=\{1,3,\dots,2d-1\}$, the sets of  even and odd integers less than $2d$, respectively.
Let $\mathcal{I}$ denote the collection of all functions $\iota: 2d\ra 2d$ such that
\begin{equation}
\{\iota(0),\iota(1)\}<\{\iota(2),\iota(3)\}<\dots<\{\iota(2d-2),\iota(2d-1)\}.
\end{equation}
Each $\iota\in \mathcal{I}$ codes two strictly increasing sequences $\iota\re D_e$ and $\iota\re D_o$, each of length $d$.
For $\vec{\gamma}\in[\kappa]^{2d}$,
$\iota(\vec{\gamma}\,)$ determines the pair of sequences of ordinals
\begin{equation}
 (\gamma_{\iota(0)},\gamma_{\iota(2)},\dots,\gamma_{\iota(2d-2))}), (\gamma_{\iota(1)},\gamma_{\iota(3)},\dots,\gamma_{\iota(2d-1)}),
\end{equation}
both of which are members of $[\kappa]^d$.
Denote these as $\iota_e(\vec\gamma\,)$ and $\iota_o(\vec\gamma\,)$, respectively.
To ease notation, let $\vec{\delta}_{\vec\al}$ denote
$\vec\delta_{p_{\vec\al}}$,
 $k_{\vec{\al}}$ denote $|\vec{\delta}_{\vec\al}|$,
and let $\ell_{\vec{\al}}$ denote  $\ell_{p_{\vec\al}}$.
Let $\lgl \delta_{\vec{\al}}(j):j<k_{\vec{\al}}\rgl$
denote the enumeration of $\vec{\delta}_{\vec\al}$
in increasing order.

Define a coloring  $g$ on $[\kappa]^{2d}$ into countably many colors as follows:
Given  $\vec\gamma\in[\kappa]^{2d}$ and
 $\iota\in\mathcal{I}$, to reduce the number of subscripts,  letting
$\vec\al$ denote $\iota_e(\vec\gamma\,)$ and $\vec\beta$ denote $\iota_o(\vec\gamma\,)$,
define
\begin{align}\label{eq.fseq}
g(\iota,\vec\gamma\,)= \,
&\lgl \iota, \varepsilon_{\vec{\al}}, k_{\vec{\al}},
p_{\vec{\al}}(x_*),
\lgl \lgl p_{\vec{\al}}(x,\delta_{\vec{\al}}(j)):j<k_{\vec{\al}}\rgl:x\in X^*\rgl,\cr
& \lgl  \lgl x,j \rgl: x\in X^*,\ j<k_{\vec{\al}}, \ \delta_{\vec{\al}}(j)=\al_x \rgl,\cr
&\lgl \lgl j,k\rgl:j<k_{\vec{\al}},\ k<k_{\vec{\beta}},\ \delta_{\vec{\al}}(j)=\delta_{\vec{\beta}}(k)\rgl\rgl.
\end{align}

Let $g(\vec{\gamma}\,)$ be the sequence $\lgl g(\iota,\vec\gamma\,):\iota\in\mathcal{I}\rgl$, where $\mathcal{I}$ is given some fixed ordering.
Since the range of $g$ is countable,
by the \Erdos-Rado Theorem 
there is a  subset $K\sse\kappa$ of cardinality $\aleph_1$
which is homogeneous for $g$.
Take $K'\sse K$ such that between each two members of $K'$ there is a member of $K$.
Take subsets $K_x\sse K'$, $x\in X^*$, such that   each $|K_x|=\aleph_0$.
 and 
$K_x<K_{x'}$ for $x<_{\mathrm{lex}} x'$ in $X^*$.
Let $\vec{K}$ denote $\prod_{x\in X^*} K_x$, and note that 
 for all $\vec\al\in \vec{K}$, $p_{\vec\al}(x_*)$ is the same level set in $E$ with 
 coding node $p_{\vec\al}(x_*)$.
 Denote this level set  by 
$Y^*$, and
 let $t^*_{x_*}$ denote the coding node   in 
$Y^*$ and  $\ell^*$ denote the length of the nodes in $Y^*$.
Note that $t^*_{x_*}= p_{\vec\al}(x_*)$ for any/all $\vec\al\in\prod_{x\in X} K_x$.

The next two lemmas  are  exactly  the same (with identical proofs) as for the forcing proof of the  Halpern-\Lauchli\ Theorem.  (See Section 3 of \cite{DobrinenJML23}.)

\begin{lem}\label{lem.HLonetypes}
There are $\varepsilon_*\in 2$, $k_*\in\mathbb{N}$,
and $ \lgl t_{x,j}: j<k_*\rgl$, $x\in X^*$,
 such that
 $\varepsilon_{\vec{\al}}=\varepsilon_*$,
$k_{\vec\al}=k_*$,   and
$\lgl p_{\vec\al}(x,\delta_{\vec\al}(j)):j<k_{\vec\al}\rgl
=
 \lgl t_{x,j}: j<k^*\rgl$,
for each $x\in X^*$,
for all $\vec\al\in \vec{K}$.
\end{lem}

\begin{lem}\label{lem.HLj=j'}
Given any $\vec\al,\vec\beta\in \vec{K}$,
if $j,j'<k^*$ and $\delta_{\vec\al}(j)=\delta_{\vec\beta}(j')$,
 then $j=j'$.
\end{lem}

For any $\vec\al\in \vec{K}$ and any $\iota\in\mathcal{I}$, there is a $\vec\gamma\in[K]^{2d}$ such that $\vec\al=\iota_o(\vec\gamma)$.
By homogeneity of $f$ and  by the first sequence in the second line of equation  (\ref{eq.fseq}), there is a strictly increasing sequence
$\lgl j_x:x\in X^*\rgl$  of members of $k^*$ such that for each $\vec\al\in \vec{K}$,
$\delta_{\vec\al}(j_x)=\al_x$.
For each $x\in X^*$, let $t^*_x$ denote $t_{x,j_x}$.
Then  for each $x\in X^*$ and each $\vec\al\in \vec{K}$,
\begin{equation}
p_{\vec\al}(x,\al_x)=p_{\vec{\al}}(x, \delta_{\vec\al}(j_x))=t_{x,j_x}=t^*_x.
\end{equation}

\begin{lem}\label{lem.HLcompat}
The set of conditions  $\{p_{\vec{\al}}:\vec{\al}\in \vec{K}\}$ is  compatible.
\end{lem}

\begin{proof}
Suppose toward a contradiction that there are $\vec\al,\vec\beta\in\vec{K}$ such that $p_{\vec\al}$ and
 $p_{\vec\beta}$ are incompatible.
Since  $g$ being homogeneous implies that 
$p_{\vec\al}(x_*)=p_{\vec\beta}(x_*)$,
the incompatibility must take place in ``$X^*$ parts" of the conditions.
By Lemma \ref{lem.HLonetypes},
for each $x\in X^*$ and $j<k^*$,
\begin{equation}
 p_{\vec{\al}}(x,\delta_{\vec{\al}}(j))
=t_{x,j}
=p_{\vec{\beta}}(x,\delta_{\vec{\beta}}(j)).
\end{equation}
Thus,
 the only way $p_{\vec\al}$ and $p_{\vec\beta}$ can be incompatible is if
there are  $x\in X^*$ and $j,j'<k^*$ such that
$\delta_{\vec\al}(j)=\delta_{\vec\beta}(j')$
but
$p_{\vec\al}(x,\delta_{\vec\al}(j))\ne p_{\vec\beta}(x,\delta_{\vec\beta}(j'))$.
Since
$p_{\vec\al}(x,\delta_{\vec\al}(j))=t_{x,j}$ and
$p_{\vec\beta}(x,\delta_{\vec\beta}(j'))= t_{x,j'}$,
this would imply
 $j\ne j'$.
But by Lemma \ref{lem.HLj=j'},
$j\ne j'$ implies that $\delta_{\vec\al}(j)\ne\delta_{\vec\beta}(j')$, a contradiction.
Therefore,
 $p_{\vec\al}$ and $p_{\vec\beta}$ must be  compatible.
\end{proof}

We now construct $T\in [d,S]$  satisfying the theorem. 
\vskip.1in

\noindent {\underline{Prepping Stage}}.
Recall that $X=X^*\cup \{x_*\}$ and let $Z=\max(r_d(S))^+\setminus X$, so that $\max(r_d(S))^+=X\cup Z$.
Recall that  $\ell^*$ denotes the length of the nodes $t^*_x$, $x\in X$.
For every  $z\in Z$, 
let  $t^*_z$ denote the leftmost extension of $z$ in $S\re \ell^*$,
and let  $W^*=\{t^*_w:w\in X\cup Z\}$.
Then $W^*$ is 
 a level set extension  of the   level set $\max(r_d(S))^+$.
 Hence, there 
   is an $a_0\in r_{n_0+1}[d,S]$  with the property that
the level set of nodes 
   $\max(r_{d+1}(a_0))$ extends the level set $W^*$.
In particular,  $r_d(a_0)=r_d(S)$ and every node in $a_0\setminus r_d(S)$ end-extends some node in $W^*$.
\vskip.1in

\noindent {\underline{Step 0}}.
For each $x\in X^*\setminus\{y\}$, let $J_x$ be a subset of $K_x$ with cardinality the number of nodes in $\max(a_0)$ extending $t^*_x$ in $U'_x$.
If $x_*\ne y$, then  let
$Z$ be the set of those $s\in\max(a_0)$ such that $t^*_y\sse s$ and $s\contains c^{\frown}1$, where $c$ is the immediate precedessor of $s$ in $a_0$. 
Let 
 $J_y$ be a subset of $K_y$ with cardinality $|Z|$.
Let $\vec{J}=\prod_{x\in X^*}J_x$.
The set of conditions $\{p_{\vec\al}:\vec\al\in\vec{J}\}$ is compatible, by Lemma \ref{lem.HLcompat}.
Let $q$ be a condition in $\bP$ such that for each $x\in X^*$,
$\{q(x,\delta):\delta\in J_x\}$ lists the nodes in $\max(a_0)$ extending $t^*_x$. 
Let
$q(x_*)$ denote   the coding node in $\max(a_0)$; note that  $q(x_*)\contains t^*_{x_*}$.
Expand the domain of $q$ to 
$\vec{\delta}_q=\bigcup_{\vec\al\in\vec{J}}\vec{\delta}_{p_{\vec{\al}}}$ 
obtaining  a condition  $q\in \bP$ such that 
$q\le p_{\vec\al}$, for each $\vec\al\in\vec{J}$.
That is, 
for each pair $(x,\gamma)$ with $x\in X^*$ and $\gamma\in\vec{\delta}_q\setminus
J_x$,
there is at least one $\vec{\al}\in\vec{J}$ and some $j'<k^*$ such that $\delta_{\vec\al}(j')=\gamma$.
For any other $\vec\beta\in\vec{J}$ for which $\gamma\in\vec{\delta}_{\vec\beta}$,
since the set $\{p_{\vec{\al}}:\vec{\al}\in\vec{J}\}$ is pairwise compatible by Lemma \ref{lem.HLcompat},
it follows
 that $p_{\vec\beta}(x,\gamma)$ must  equal $p_{\vec{\al}}(x,\gamma)$, which is exactly $t^*_{x,j'}$.
Let $q(x,\gamma)$ be the leftmost extension
 of $t_{x,j'}^*$ in $S$.
Thus, $q(x,\gamma)$ is defined for each pair $(x,\gamma)\in d\times \vec{\delta}_q$.
Define
\begin{equation}
q= q(x_*)\cup \{\lgl (x,\delta),q(x,\delta)\rgl: x\in X^*,\  \delta\in \vec{\delta}_q\}.
\end{equation}

By the  Truth Lemma, there is some $r\le q$ in $\bP$ such that  $\varepsilon(r)=\varepsilon_*$.
Take $b_0$ to be a member of $r_{n_0+1}[d,S]$ such that 
\begin{enumerate}
    \item
Each node in $\max(b_0)$ 
end-extends some node in $\max(a_0)$;
\item
For each $x\in X^*$ and each $\gamma\in J_x$, 
$r(x,\gamma)$ is the node in $\max(b_0)$ extending the node $q(x,\gamma)$ in $\max(a_0)$;
\item
The node in $\max(b_0)$ extending $q(x_*)$ is $r(x_*)$.
\end{enumerate}
\ \vskip.1in

\noindent {\underline{Induction Hypothesis}}. 
Given $i\ge 1$ suppose that we have built $b_0,\dots, b_{i-1}$ so that 
$b_0\in r_{n_0+1}[d,S]$,
for each $0<j<i$, $b_j\in r_{n_j+1}[b_{j-1},S]$, 
and 
for each $j<i$, we have 
 $c^{b_{i-1}}_{n_j}$ (which equals $c^{b_{j}}_{n_j}$) is a coding node in $M$ which extends $x_*$  and
 $f$ has constant value 
 $\varepsilon_*$ on $ (E\cap b_{i-1})\otimes \bigotimes_{x\in X_0} (U_x\cap b_{i-1})$.
\vskip.1in

\noindent {\underline{Induction Step $i$}}. 
This is  essentially the same as Step 0.
Take  $a_i$ to  be a member of $r_{n_i+1}[b_{i-1},S]$.
(Note that for each $x\in X$,  all nodes in $a_i$ extending $x$  also extend $t^*_x$.)
For each $x\in X^*\setminus\{y\}$, let $J_x$ be a subset of $K_x$ with cardinality the number of nodes in $\max(a_i)$ extending $t^*_x$ in $U'_x$.
If $x_*\ne y$, then  let
$Z$ be the set of those $s\in\max(a_i)$ such that $t^*_y\sse s$ and $s\contains c^{\frown}1$, where $c$ is the immediate precedessor of $s$ in $a_i$. 
Let 
 $J_y$ be a subset of $K_y$ with cardinality $|Z|$.
Let $\vec{J}=\prod_{x\in X^*}J_x$.
The set of conditions $\{p_{\vec\al}:\vec\al\in\vec{J}\}$ is compatible, by Lemma \ref{lem.HLcompat}.
Let $q$ be a condition in $\bP$ such that for each $x\in X^*$,
$\{q(x,\delta):\delta\in J_x\}$ lists the nodes in $\max(a_i)$ extending $t^*_x$. 
Let
$q(x_*)$ denote   the coding node in $\max(a_i)$; note that  $q(x_*)\contains t^*_{x_*}$.
Expand the domain of $q$ to 
$\vec{\delta}_q=\bigcup_{\vec\al\in\vec{J}}\vec{\delta}_{p_{\vec{\al}}}$ 
obtaining  a condition  $q\in \bP$ such that 
$q\le p_{\vec\al}$, for each $\vec\al\in\vec{J}$.
Then by
 the Truth Lemma,  there is some 
$r\le q$ in $\bP$ such that  
$f(\{r(x_*)\}\cup\{ r(x,\delta_x):x\in X^*\})=\varepsilon_*$ for each sequence $\langle \delta_x:x\in X^*\rangle\in\prod_{x\in X^*}J_x$.

Take $b_i$  in $r_{n_i+1}[b_{i-1},S]$ such that 
\begin{enumerate}
    \item
Each node in $\max(b_i)$ 
end-extends some node in $\max(a_i)$;
\item
For each $x\in X^*$ and each $\gamma\in J_x$, 
$r(x,\gamma)$ is the node in $\max(b_i)$ extending the node $q(x,\gamma)$ in $\max(a_i)$;
\item
The node in $\max(b_i)$ extending $q(x_*)$ is $r(x_*)$.
\end{enumerate}
\ \vskip.1in

At the end of this induction, let $T=\bigcup_{i<\om} b_i$.  Then $T\in [d,S]$ and  for each  coding node $c\in E\cap T$,
all vectors in $\{c\}\times \prod_{x\in X^*} (T\cap U_x)\re |c|$
have the same color.
By the End-Homogeneity,  the coloring is homogeneous on  
$\bigotimes_{x\in X,x_*} V_x$,
where $V_x:= U_x\cap T$ for $x\in X$.
Hence,
 $T$ satisfies the conclusion of the theorem.
\end{proof}

As a corollary, we recover  Proposition 7 in  \cite{CEW_EUROCOMB25} that $\Psi$ is indivisible. 

\begin{cor}
For each coloring of the nodes in $\Psi$ into finitely many colors, 
there is a subcopy of $\Psi$ in which all nodes have the same color.
\end{cor}

\begin{proof}
This can be seen two ways:
Apply Theorem \ref{thm.HLVariant} to the special case  $d=0, X_0=\emptyset, X_2=\{y\}$.
For a proof that does not use forcing techniques,  note that the collection of coding nodes that are minimal amongst those with its $\theta$-value are dense in the subtree obtained by 
 Theorem \ref{lem.Stepone} applied to the special case 
$d=0, X_0=\emptyset, X_2=\{y\}$; then
apply the fact that the rationals are indivisible.
\end{proof}

\subsection{A Milliken Theorem for Diaries}\label{subsec.DiaryMilliken}

This subsection is focused on proving a Ramsey theorem for diaries for finite chains in $\Psi$, Theorem \ref{thm.FBRD}.
For a  diary $\Delta$ for a finite chain in $\Psi$ and $m\le \height(\Delta)$,
$\Delta|_m$ denotes $\Delta\cap 3^{<m}$, the first $m$ many levels of $\Delta$.
In particular, $\Delta|_0$ is the emptyset, $\Delta|_1$ is a single critical node of type (a)--(c), and $\Delta|_{\height(\Delta)}=\Delta$.
For $\ell<\height(\Delta)$, let 
$\Delta(\ell)$ denote $\Delta\cap 3^{\ell}$, the set of nodes in $\Delta$ of length $\ell$. 
A function $f$ is a member of $\Emb(\Delta|_m,\bS)$
iff $f$ is the truncation to $\Delta|_m$ of some embedding in $\Emb(\Delta,\bS)$.
For such an $f$ and any $\ell<m$, let   $\tilde{f}(\ell)$ denote the 
length of the nodes in   $f[\Delta(\ell)]$.

We start this subsection  with the following  corollary of Theorem \ref{thm.HLVariant}, phrased in terms of  diaries and embeddings.

\begin{cor}\label{cor.End-homog}
If $\Delta$ is a diary for a finite chain in $\Psi$, $m<\height(\Delta)$, and $f\in\Emb(\Delta|_m,\bS)$ is fixed, then for any coloring $\chi$ of $G:=\{g\in\Emb(\Delta|_{m+1},\bS)  :g|_m=f\}$ into finitely many colors, there is a $\phi\in\Emb(\bS,\bS)$ which is the identity on the first $\tilde{f}(m-1)$ levels and with $\{\phi\circ g:g\in G\}$ monochromatic for $\chi$.
\end{cor}

\begin{proof}
Since indivisibility is already established, we may assume that $\height(\Delta)\ge 2$.
Let
$d_m$ denote the critical node in $\Delta(m)$.
Fix the following notation to align with Assumption \ref{assum.prodtrees}.

\underline{Case 1}.
 $m=0$.  Then $\Delta|_0=\emptyset$,  and  $\Delta|_1=\{d_0\}$ is a singleton critical node of type (a) or  (c).
In this case, $x_*=y=c_0$,   $X=X_2=\{y\}$, and   $U_y=\bS$.

\underline{Case 2}.
$m>0$.  Let $\tilde{m}=\tilde{f}(m-1)+1$ and let
$X=g[\Delta(m)]\re \tilde{m}$ for any/all $g\in G$.
 Let   $x_*$ denote the node in $X$ such that $x_*\sse g(d_m)$  for any/all    $g\in G$.
Let $d_{m-1}$ denote the critical node in $\Delta(m-1)$.  If $d_{m-1}$ is a critical node of type (c) or (d) then 
$y = g(s)\re \tilde{m}$, where $s$ is  the immediate successor of $d_{m-1}$ in $\Delta$, and  $X_2=\{y\}$.
If $d_{m-1}$ is a critical node of type (a) or (b) then 
$X_2=\emptyset$, since  all nodes
 in $\Delta(m)$ extend a node in  $\Delta(m-1)$  by  either a $0$ or a $2$.

Let $k$ be least such that $f[\Delta|_m]\sse r_k(\mathbb{S})$.
Apply Theorem
 \ref{thm.HLVariant}  to obtain an $S\in [k,\bS]$ so that  all $g\in G$ with $\ran(g)\sse S$ have the same $\chi$-color.
Let $\phi$ denote the member of $\Emb(\bS,\bS)$ such that $\phi[\bS]=S$.
Then the corollary holds. 
\end{proof}

In what follows,  for $S\in\mathcal{S}$ and $k<\om$,  
let $\Emb(\Delta|_m,r_k(S))$ denote the set of those $f\in\Emb(\Delta|_m,S)$ such that $f[\Delta(m-1)]$ are maximal nodes in $r_k(S)$.
Given an 
$f\in\Emb(\Delta|_m,\bS)$,
let
$G_f$ denote $\{g\in\Emb(\Delta|_{m+1},\bS)  :g|_m=f\}$.

\begin{lem}[End-homogeneity]\label{lem.End-homog}
Suppose $\Delta$ is a diary for a finite chain in $\Psi$ and  $m<\height(\Delta)$.  Then for any coloring $\chi$ of  $\Emb(\Delta|_{m+1},\bS)$, there is a 
 $\phi\in\Emb(\bS,\bS)$ such that 
for each $f\in\Emb(\Delta|_m,\bS)$, $\chi$ is monochromatic on 
$\{\phi\circ g:g\in G_f\}$, where 
$G_f:=\{g\in\Emb(\Delta|_{m+1},\bS)  :g|_m=f\}$.
\end{lem}

\begin{proof}
This is a standard fusion argument.
Let $\langle k_i: i<\om\rangle $ enumerate those $k<\om$ for which 
$\Emb(\Delta|_m,r_k(\bS))$ is non-empty.

For the basis of the fusion argument, 
list the members of $\Emb(\Delta|_m, r_{k_0}(\bS))$ as $f_n$, $n\le n_0$.
These satisfy $\tilde{f}_n(m-1)=k_0-1$.
Apply Corollary \ref{cor.End-homog} to obtain a $\phi^0_0\in\Emb(\bS,\bS)$ such that 
$\phi^0_0$ is the identity on $r_{k_0}(\bS)$ and 
$\chi$ is monochromatic on $\{\phi^0_0\circ g:g\in G_{f_0}\}$.
Next, apply Corollary \ref{cor.End-homog} to obtain a $\phi^1_0\in\Emb(\bS,\bS)$ such that
$\phi^1_0$ is the identity on $r_{k_0}(\bS)$ and 
 $\chi$ is monochromatic on $\{\phi^1_0\circ \phi^0_0\circ g:g\in G_{f_1}\}$.
Continuing in this  manner applying  Corollary \ref{cor.End-homog},
at the end of the 
induction  we obtain   a 
 $\phi^{n_0}_0\in\Emb(\bS,\bS)$ such that 
$\phi^{n_0}_0$ is the identity on $r_{k_0}(\bS)$ and 
$\chi$ is monochromatic on $\{\phi^{n_0}_0\circ\cdots\circ\phi^1_0\circ\phi^0_0\circ g:g\in G_{f_{n_0}}\}$.
Let $\phi_0=\phi^{n_0}_0\circ\cdots\circ\phi^1_0\circ\phi^0_0$, and let $S_0=\phi_0[\bS]$.
Then  $\phi_0$ is  the identity on $r_{k_0}(\bS)$, so   $S_0\in [r_{k_0},\bS]$.
Further, for each $n\le n_0$, $\chi$ is monochromatic on 
$\{ \phi_0\circ g:g\in G_{f_n}\}$.

 Let $S_{-1}$ denote $\bS$ and $\phi_{-1}$ denote the identity on $\bS$, and let $j\ge 0$.
Suppose 
for each $0\le i\le j$,
we have found $\phi_i\in\Emb(\bS,\bS)$ so that, letting $S_i=\phi_i[\bS]$, 
 the following hold:
\begin{enumerate}
\item
$\phi_i$ is the  identity on $r_{k_i}(S_{i-1})$;
\item
$\phi_i=\phi'\circ \phi_{i-1}$ for some $\phi'\in\Emb(\bS,\bS)$;

\item
For each 
 $f\in \Emb(\Delta|_m, r_{k_i}(S_{i-1}))$,
$\chi$ is monochromatic on $\{\phi_i\circ g:g\in G_f\}$.
\end{enumerate}
It follows that for each $0\le i\le j$ and each 
$f\in \Emb(\Delta|_m, r_{k_i}(S_j))$,
$\chi$ is monochromatic on $\{\phi_j\circ g:g\in G_f\}$.

For the induction step,
as in the base case  apply
 Corollary \ref{cor.End-homog} successively to each  member of $\Emb(\Delta|_m, r_{k_{j+1}}(S_{j}))$ to obtain a $\phi_{j+1}$
which is the identity on $r_{k_{j+1}}(S_{j})$
so that  for each $f\in\Emb(\Delta|_m, r_{k_{j+1}}(S_{j}))$, $\chi$ is monochromatic on 
$\{ \phi_{j+1}\circ g:g\in G_f\}$.

To finish, let $T=\bigcup_{i<\om}r_{k_i}(S_{i-1})$, and let $\phi\in \Emb(\bS,\bS)$ be the embedding such that for each $i<\om$, $\phi\re r_{k_i}(\bS)= \phi_i\re r_{k_i}(\bS)$.
Then  $\phi[\bS]=T$ and $\phi$ 
satisfies the lemma.
\end{proof}

The final theorem of this Section shows that initial segments of diaries have the Ramsey property.
It then follows that each diary has the Ramsey property.

\begin{thm}\label{thm.FBRD}
Let $\Delta$ be a diary for a finite chain in $\Psi$, and fix  $0\le m<\height(\Delta)$.
For any coloring $\chi$ of $\Emb(\Delta|_{m+1},\bS)$  into finitely many colors, there is a $\phi\in\Emb(\bS,\bS)$ 
with $\{\phi\circ g:g\in   \Emb(\Delta|_{m+1},\bS)  \}$ monochromatic for $\chi$.
In particular,  if $\chi$ is a finite coloring of $\Emb(\Delta,\bS)$, then there is a $\phi\in\Emb(\bS,\bS)$ with $\chi$ constant on $\phi\circ \Emb(\Delta,\bS)$.
\end{thm}

\begin{proof}
The proof is by induction on $m$.
If $m=0$, then 
Lemma \ref{lem.End-homog} yields a $\phi\in\Emb(\bS,\bS)$ such that 
$\chi$ is monochromatic on $\{\phi\circ g :g\in \Emb(\Delta|_1,\bS)\}$.

Now assume that $m\ge 1$ and the theorem holds for all $m'<m$.
By Lemma \ref{lem.End-homog}, there is a $\phi_0\in \Emb(\bS,\bS)$ so that $\chi$ is monochromatic on 
$\{\phi_0\circ g:g\in G_f\}$, for each $f\in \Emb(\Delta|_m,\bS)$.
This induces a coloring $\chi'$ on 
$\Emb(\Delta|_m,\bS)$, where for 
$f\in \Emb(\Delta|_m,\bS)$,
$\chi'(f)$ is the  constant color of $\chi$ on 
$\{\phi_0\circ g:g\in G_f\}$.
By the induction hypothesis, there is a $\phi_1\in \Emb(\bS,\bS)$ such that 
$\chi'$ is constant on 
$\phi\circ \Emb(\Delta|_m,\bS)$, where $\phi=\phi_1\circ \phi_0$ for some $\phi_1\in\Emb(\bS,\bS)$.
Then $\chi$ is constant on $\phi\circ\Emb(\Delta|_{m+1},\bS)$.

The last sentence follows from the case that $m+1=\height(\Delta)$.
\end{proof}

\section{Proof of Main Theorem}\label{sec.mainresults}

This final section of the paper contains of a proof of the Main Theorem (Subsection \ref{section.aa}), the seven diaries for chains in $\Psi$  of length two (Subsection \ref{subsec.ex7}), and  a discussion of further directions (Subsection \ref{subsec.furtherdir}).

\subsection{Proof of the Main Theorem}\label{section.aa}

Recall Definition \ref{defn.aac}:
An almost antichain is a subset $A$ of coding nodes in $\bS$ such that  (a)
each pair of coding nodes $c\ne d$ in $A$ is either 
incomparable,
or else  $c^{\frown}1\sse d$ or $d^{\frown}1\sse c$; and (b) if $s$ is the meet of two incomparable nodes in $A$, then each node in $A$ extending $s$ extends either  $s^{\frown}0$ or $s^{\frown}2$.
Recall that 
$\Psi^*|_A$ denotes  the restriction of $\Psi^*$ to the set of nodes $\{p_i\in\Psi^*: i\in I(A)\}$, where 
$I(A)$ is the set of those indices $i<\om$ for which the coding node $c_i\in\bS$ is a coding node in $A$.
We will say that an almost antichain $A$ {\em encodes a copy of} $\Psi^*$ iff $\Psi^*|_A$
is isomorphic to $\Psi^*$.

In Lemma \ref{lem.existsaac},
we show that in each $S\in\mathcal{S}$ contains an almost antichain $A\sse S$ encoding a copy of $\Psi^*$.
These are 
 minimalistic in the sense that  each chain in $\Psi^*|_A$ 
uniquely induces a 
diary-shaped subset of  $A$, as will be shown in Theorem \ref{thm.Asubsetsdiaries}.
The Main Theorem will then follow immediately from Theorem \ref{thm.FBRD},  Lemma \ref{lem.existsaac}, and 
Theorem \ref{thm.Asubsetsdiaries}.
We remark that, in contrast to all other known finite big Ramsey degree results,  
an actual antichain  of coding nodes in $\bS$ cannot represent a subcopy of $\Psi^*$, since an antichain $B\sse \bS$ 
omits  some coding nodes representing meets of two nodes on different rays in $\Psi^*|_B$ and hence, $\Psi^*|_B$ will not be closed under meets.

\begin{lem}\label{lem.existsaac}
Given any $S\in\mathcal{S}(\bS)$, there is an almost antichain 
 $A\sse S$ encoding a copy of $\Psi^*$.
\end{lem}

\begin{proof}
Recall 
that $S\in \mathcal{S}(\bS)$ implies that $S$
is isomorphic to $\bS$ and hence, represents a subcopy of $\Psi^*$ in the same enumeration as $\bS$.
The following notation will be useful.
For $n<\om$, let $\bS(n)$ denote $\{s\in \bS:|s|=|c_n|\}$,  where $c_n$ denotes the $n$-th coding node in $\bS$.
Recall  that  $\bS(n)$ has $2n+1$ nodes.
For $n<\om$, let $\{s_{n,i}:i<2n+1\}$ enumerate the nodes in  $\bS(n)$  in lexicographic order.

Recall that the proof of Observation \ref{obs.2.9} only used (c), (e), and (f) of Definition \ref{defn.ctiso}.
We will 
construct nested increasing  almost antichains  $A_0\sse A_1\sse \dots\sse S$, where
$A_n$ has coding nodes $\{c^A_i:i\le n\}$ and
$A_{n+1}$ end-extends $A_n$, so that
 $A:=\bigcup_{n<\om}A_n$ is an almost antichain  and there is an order-preserving bijection between  $\bS(n+1)$ and  the set of immediate successors in $A$ of the nodes 
of length $|c_n^A|$
 satisfying  (c), (e), and (f) 
of Definition \ref{defn.ctiso}.
The  same proof as that of  Observation \ref{obs.2.9}  will then yield lemma.

Construction of $A_0$:
Let $u_0=c^S_0$;
this is the root of $A$ and also a splitting (but not coding) node in $A$.
Let $v_0$ be the least coding node  in $S$  extending 
${u_0}^{\frown}0$, noting  that  $\theta(v_0)=\theta(u_0)$.
This $v_0$ is  the next 
 splitting (but not coding) node in $A$.
Let $c^A_0$ denote the least coding node in $S$  extending ${v_0}^{\frown} 2$, noting that $\theta(c^{A}_0)=\theta(u_0)$; $c^A_0$ is the least coding node in $A$.
Let $a_{0,0}$ denote the leftmost extension  of ${v_0}^{\frown} 0$
 in $S$ with $|a_{0,0}|=|c^A_0|$, and let $a_{0,2}$ denote the leftmost extension of ${u_0}^{\frown} 2$ with 
$|a_{0,2}|=|c^A_0|$.
Let $a_{0,1}$ denote $c^A_0$ and
$\theta_0'=\theta(c_0^A)$,
and note that
\begin{equation}\label{eq.16}
\theta(a_{0,0})=\theta(a_{0,1})=\theta(c^A_0)=\theta(a_{0,2})=
\theta(u_0)=\theta(v_0)=\theta'_0.
\end{equation}
The $0$-the coding level of $A$ is 
\begin{equation}
\CL_A(0):=\{a_{0,0}, a_{0,1}, a_{0,2}\},
\end{equation}
and 
\begin{equation}
A_0:=\{s\re \ell: s\in \CL_A(0)\wedge \ell\in\{|u_0|,|v_0|,|c^A_0|\}\}
\end{equation}
is easily seen to be  an almost antichain
with critical nodes $\Crit(A_0)=\{u_0,v_0,c^A_0\}$.

Define $\IS_A(0)$, the set of immediate successors in $\widehat{A}$ of $\CL_A(0)$, as follows: 
 ${a_{0,1}}^+={a_{0,1}}^\frown 1$, and ${a_{0,0}}^+,{a_{0,2}}^+$  are the immediate successors of $a_{0,0},a_{0,2}$ in $\widehat{S}$, respectively.
Then $\IS_A(0)=\{{a_{0,j}}^+:j<3\}$.
Note that ${a_{0,0}}^+= {a_{0,0}}^{\frown}0$, since every coding node $c$ in $S$ extending $ {a_{0,0}}^+$  represents a point in $\Psi^*$ that is  $\prec$-above and $<_{\mathrm{lex}}$ left of the point in $\Psi^*$ represented by $c^A_0$.
Similarly, 
 ${a_{0,2}}^+={a_{0,2}}^{\frown}2$,
since every coding node  $c$ in $S$ extending $ {a_{0,2}}^+$  with the same $\theta$-value 
represents a point in $\Psi^*$ that is  
 $\prec$-below   the point in $\Psi^*$ represented by  $c^A_0$. 
Hence,  for each $j<3$, $a_{0,j}^+(|c^A_0|)=s_{0,j}(|c_0|)$.
Let $\varphi_0:\bS(1)\ra \IS_A(0)$ be the lexicographic order-preserving isomorphism. 
Then  the Induction Hypothesis,  below, is satisfied for $n=0$.

Let $n$ be fixed.
Suppose  we have constructed  end-extending almost antichains $A_0\sse \dots\sse A_n\sse S$ so that  for each $k\le n$,
$\Crit(A_k)=\{u_i,v_i,c^A_i:i\le k\}$ 
is  the set of critical nodes in $A_k$, 
$\CL_A(k)=\{a_{k,j}:j<2(k+1)+1\}$ is the lexicographic-ordered enumeration
of  the nodes in $A_k$ with length $|c^A_k|$, 
and $\IS_A(k)=\{{a_{k,j}}^+:j<2(k+1)+1\}$ is the set of 
immediate successors in $\widehat{A}$ of the nodes in $\CL_A(k)$. 
Let $\varphi_i$ be the lexicographic order-preserving bijection from 
 $\bS(i+1)$ to 
$\IS_A(i)$,  and let $\tilde{\varphi}_k=\bigcup_{i\le k}\varphi_i$.
Assume  the following Induction Hypothesis  holds for   each $k\le n$:
\begin{enumerate}
\item[(1)]
$A_k:=\{s\re |w|: s\in \CL_A(n)\wedge
w\in \Crit(A_k)\}$
 is  an almost antichain, with $a_{k,j}=c^A_k$ iff $s_{k+1,j}={c_k}^{\frown}1$, for $j<2(k+1)+1$.

\item[(2)]

If $a_{k,j}\ne c^A_k$, then every node  $t$ on the leftmost branch (and on the rightmost branch) of $S$ extending ${a_{k,j}}^+$ satisfies  $\theta(t)=\theta({a_{k,j}})=\theta({a_{k,j}}^+)$.

\item[(3)]
$\Psi^*|_{A_k}$ is isomorphic to $T_k$.

\item[$(c')$]

${a_{k,j}}^+(|c^A_k|)=s_{k+1,j}(|c_k|)$, 
for   each  $j<2(k+1)+1$.

\item[$(e')$]
For each $a,b\in \bigcup_{i\le k}\IS_A(i)$,
$\theta(a)\le\theta(b)$ iff $\theta(\tilde{\varphi}_k(a))\le\theta(\tilde{\varphi}_k(b))$.

\item[$(f')$]
For each $a,b\in \bigcup_{i\le k}\IS_A(i)$,
$|a|\le |b|$ iff $|\tilde{\varphi}_k(a)|\le\ |\tilde{\varphi}_k(b)|$.
\end{enumerate}

For the Inductive Step, let  $m=n+1$, and let $\ell<2m+1$ be the index such that $c_n^A=a_{n,\ell}$.
Take $c_*$ to be the least coding node in $S$ extending ${a_{n,\ell}}^+:={a_{n,\ell}}^{\frown}1$
 such that $\theta(c_*)>\theta(a)$ for all $a\in\CL_A(n)$, and  let  $b_{\ell}=  {c_*}^{\frown}0$.
For each $j<2m+1$ not equal to $\ell$, let ${b_{j}}$ be the leftmost extension of ${a_{n,j}}$ in $\widehat{S}$ of length $|{b_{\ell}}|$.
Since (2) holds for $n$, 
 $j\ne\ell$ implies
all extensions of  ${a_{n,j}}^+$ along the leftmost (and rightmost) branches of $S$ have the same $\theta$-value as ${a_{n,j}}^+$.
These steps set us up so that the rest of the construction will ensure (2) and ($e'$)  hold for $m$.

Let $i_*<2m+1$ be the index such that $s_{m,i_*}=c_{m}$, 
and recall that  $i_*\ne\ell$ and 
$c_m\not\contains {c_n}^{\frown}1$.
By  ($c'$)  of the Induction Hypothesis,
${a_{n,i_*}}^+(|c^A_n|)= s_{m,i_*}(|c_n|)\ne 1$, so 
$\theta({a_{n,i_*}}^+)=\theta({a_{n,i_*}})$.
Hence, also $\theta({b_{i_*}})=\theta({a_{n,i_*}})$.
Let $u_m$ be the least coding node in $A$ extending 
${b_{i_*}}$; let
$v_m$ be the least coding node in $A$ extending ${u_m}^{\frown}0$
and $c^A_m$ be the least coding node in $A$ extending 
${v_m}^{\frown}2$.
Note that $u_m$ and $v_m$ are both  on the leftmost branch extending $b_{i_*}$, so by (2) of the Induction Hypothesis,
 $\theta'_m:=\theta(c^A_m)
=\theta(v_m)=\theta(u_m)=\theta(b_{i_*})$.
Let $a_{m,i_*+1}=c^A_m$;
let $a_{m,i_*}$   be the leftmost  extension of ${v_m}^{\frown}0$ in $A$ of length $|c^A_m|$; and let 
$a_{m,i_*+2}$   be the leftmost  extension of ${u_m}^{\frown}2$ in $A$ of length $|c^A_m|$.
Then 
\begin{equation}
\theta(a_{m,i_*})=\theta(a_{m,i_*+1})=
\theta(c^A_m)
=\theta(a_{m,i_*+2})=
\theta(v_m)=\theta(u_m) =\theta({a_{n,i_*}}^+)=\theta'_m,
\end{equation}
and each leftmost (or rightmost) extension of $u_m$  or $v_m$ retains the same $\theta$-value.  
This is enough to ensure that (2) of the Induction Hypothesis will hold for $m$.

Let  $B=\{{b_{j}}:j<2m+1,\ j\ne i_*\}$.
For each $b\in B$,
let $b'$ be the leftmost extension of $b$ in $S$ of length $|c^A_m|$, and note that $\theta(b')=\theta(b)$.
Define 
\begin{equation}
\CL_A(m):=
\{a_{m,i_*}, a_{m,i_*+1},a_{m,i_*+2}\}\cup\{b':b\in B\}.
\end{equation}
For $j<i_*$, let $a_{m,j}= (b_{n,j})'$.
For $j>i_*$, let $a_{m,j+2}=(b_{n,j})'$.
Then    
$\{a_{m,j}:j<2(m+1)+1\}$ 
enumerates 
 $\CL_A(m)$ in lexicographic order.
Let ${a_{m,i_*+1}}^+= {a_{m,i_*+1}}^{\frown}1$.
For each $j<2(m+1)+1$ with $j\ne i_*+1$,
 let 
${a_{m,j}}^+$ be the immediate successor of 
$a_{m,j}$ in $\widehat{S}$ of length $|c^A_m|+1$.  
This constructs  the set of immediate successors, $\IS_A(m):=\{a_{m,j}^+:j<2(m+1)+1\}$, of the nodes in $\CL_A(m)$.
Define 
\begin{equation}
A_m:= A_n\cup \{s\re |w|: s\in \CL_A(m)\wedge
w\in \{u_m,v_m,c^A_m\}.
\end{equation}
Let $\varphi_m:\bS(m+1)\ra\IS_A(m)$ be the order-preserving bijection, and let $\tilde{\varphi}_m=\tilde{\varphi}_n\cup \varphi_m$.

We  check that the Induction Hypothesis  holds for $m$.
By construction, $A_m$ is an almost antichain
and $a_{m,i_*+1}=c^A_m$ while $s_{m+1,i_*+1}= {c_m}^{\frown}1$,
 so (1) holds.
Also by construction, (2), ($e'$), and ($f'$) hold for $m$.
By ($c'$) for $n$, 
($e'$) and ($f'$)  for $m$,  and the proof of Observation \ref{obs.2.9}.
$\Psi^*|_{A_m}$ is isomorphic to $T_m$, 
It remains to check that  ($c'$) holds for $m$.


For  $j\in \{i_*,i_*+1,i_*+2\}$, ${a_{m,j}}^+(|c^A_m|)=s_{m+1,j}(|c_m|)$.
Now fix $j<2(m+1)+1$ with $j\not\in \{i_*,i_*+1,i_*+2\}$, and 
note that $s_{m+1,j}(|c_m|)\ne 1$.
By ($c'$) for $n$ and the rest of the Induction Hypothesis for $m$, we have that 
$c^A_m(|c^A_k|)=c_m(|c_k|)$  and
${a_{m,j}}^+(|c^A_k|)=s_{m+1,j}(|c_k|)$,
for all $k<m$.
If  $p_m$ is $\prec$-maximal in $T_m$, then  
$p_m^A$  is maximal in 
$\Psi^*|_{A_m}$ so 
$a_{m,j}^+(|c_m^A|)=2 = s_{m+1,j}(|c_m|)$; hence,  ($c'$) holds for $m$.

Now suppose that   $p_m$ is not $\prec$-maximal in $T_m$, and let 
$i<m$ be such that $p_i$ is the $\prec$-immediate successor of $p_m$ in $T_m$.
Then  
\begin{equation}
s_{m+1,j}(|c_m|)=0\Llra s_{m+1,j}(|c_k|)=0\Llra a_{m,j}(|c^A_k|)=0,
\end{equation}
where the last bi-implication follows from the  Induction Hypothesis. 
Since  $\Psi^*|_{A_m}$ is isomorphic to $T_m$, we have 
$a_{m,j}(|c^A_k|)=0$ iff ${a_{m,j}}^+(|c^A_m|)=0$.
Similarly, $s_{m+1,j}(|c_m|)=2$ iff ${a_{m,j}}^+(|c^A_m|)=2$.
Hence,  ($c'$) holds for $m$.

Let $A=\bigcup_{n<\om}A_n$.  Then $A$ is an almost antichain and  encodes a copy of $\Psi^*$, since each $A_n$ encodes a tree isomorphic to $T_n$.
\end{proof}

The following observation will be used in the proof of Theorem \ref{thm.Asubsetsdiaries}.

\begin{observation}
Let $C$ be a finite  almost antichain  in $\bS$ such that 
$\Psi^*|_C$ is a  chain in $\Psi^*$.
 The meet-closure  of $C$ in $\bS$,  denoted by $M_C$, has the following properties:
\begin{enumerate}
\item
No two nodes in $M_C$ have the same length. (This follows from  the fact that every node in $M_C$ is a coding node in $\bS$.)
\item
If $d$ is a node in $M_C\setminus C$, then $d$ is the meet of two nodes in $C$ and 
the immediate successors of $d$ in $\widehat{C}$ are 
$d^{\frown}0$ and $d^{\frown}2$.
(This is because $A$ is an almost antichain.)
\item
Only nodes in the leftmost path in $M_C$ can possibly have different  $\theta$-values from their immediate predecessors in $M_C$.
(This is Proposition \ref{prop.no01splits}.)
\end{enumerate}
\end{observation}

\begin{thm}\label{thm.Asubsetsdiaries}
Let $A\sse \bS$ be an almost antichain representing a copy of $\Psi^*$. Suppose  $\mathbf{C}$ is a finite chain in $\Psi^*|_A$, and 
let  $C=\{c_i\in\bS:p_i\in\mathbf{C}\}$.
Then there is a unique diary $\Delta$ and unique diary embedding $e:\Delta\ra\widehat{C}$.
In particular, $e[\Delta]$ is the unique diary-shaped subset of $\widehat{C}$.
\end{thm}

\begin{proof}
We will show that $C$ induces  a unique diary-shaped subset $D$ of $A$.
 This will immediately  induce the diary $\Delta$ and diary-embedding $e$.

First notice that $C$ is an almost-antichain, since $C\sse A$.
Let  $\langle c^A_i:i<p\rangle$ list the  nodes in $C$ in increasing order of length.
Let $M_C$ denote the meet-closure of $C$,
and 
let $\langle t_k:k<m\rangle$ enumerate the nodes in $M_C$ in increasing order of length.
Note that   $t_{m-1}=c^A_{p-1}$.
The critical nodes in $D$ will consist of the nodes in 
$M_C$ along with 
 all first instances of $\theta$-value changes that occur in the leftmost branch of $\widehat{C}$ between levels in $M_C$.

Let $i_*<p$ be the index of the lexicographic leftmost  coding node in $C$.
Any changes in $\theta$-values in $\widehat{C}$ must take place on  initial segments of $c_{i_*}$.
Let $K$ denote the set of those $k<i_*$ such that 
$\theta(c_{i_*}\re |t_k|)\ne \theta(c_{i_*}\re |t_{k+1}|)$ and 
$c_{i_*}\re \ell_k$ is not in  $C$.
Let  $b_k$
denote  $c_{i_*}\re\ell$, 
where $\ell$ is least in  the interval $(|t_k|,|t_{k+1}|)$ such that  $c_{i_*}(\ell)=1$.
Then $b_k$ is a critical node of type (d).
Let $B=\{b_k:k\in K\}$, 
$L=\{|t|:t\in M_C\cup  B\}$, and 
\begin{equation}
D=\{c\re \ell:c\in C\wedge \ell\in L\}.
\end{equation}
Then $\Crit(D)=M_C\cup B$.
Let $\langle u_j:j<n\rangle$ be the enumeration of $\Crit(D)$ in order of increasing length.
The diary
 $\Delta$ will  be of height $n$ and  determined by     $\Crit(D)$.
As in previous sections,  $\Delta|_j$  denotes $\Delta\cap 3^{\le j}$.  Here, we will let $\Delta_j$ denote $\Delta\cap 3^j$.
We now precisely define $\Delta$ and $e$.

Let $d_0$ denote  the empty sequence.
If $u_0$ is a coding node in $D$, then  $d_0$ is declared to be a coding node and we let ${d_0}^{\frown}1$ be its immediate successor.
If $u_0$ is a splitting node in $D$, then 
$d_0$ is declared to be a splitting node with 
 immediate successors  ${d_0}^{\frown}0$ and ${d_0}^{\frown}2$.
This determines the type of the critical node $d_0$ and constructs the nodes in $\Delta|_1$.
Define $e(d_0)= u_0$.
Notice that  the lexicographic order-preserving bijection $\lambda$ from $\Delta_1$ to the set of immediate successors of $u_0$ in $\widehat{D}$ satisfies $\lambda(d)(|u_0|)=d(0)$ for each $d\in \Delta_1$.

Suppose  $j<n-1$ and we are given the nodes in $\Delta|_{j+1}$ along with the types of the critical nodes $d_0,\dots, d_j$ and a diary embedding
$e:\Delta|_j\ra \{u\in D:|u|\le |u_j|\}$.
Further assume that there is a lexicographic order-preserving bijection $\lambda:\{d\in \Delta_{j+1}\}\ra \{s\in \widehat{D}:|s|=|u_j|+1\}$, and this $\lambda$ preserves passing numbers, meaning that $d(j)=\lambda(d)(|u_j|)$.
Let $d_{j+1}$ be the node in $\Delta_{j+1}$ such that $\lambda(d_{j+1})=u_{j+1}\re(|u_j|+1)$.
Declare $d_{j+1}$ to be a critical node of type (a)  iff $u_{j+1}$ is a splitting node in $D$ (i.e., $u_{j+1}\in M_C\setminus C$);
type (b) iff $u_{j+1}$ is a terminal coding node in $C$;
type (c) iff $u_{j+1}$ is a non-terminal coding node in $C$;
type (d) iff $u_{j+1}\in B$. 
This defines $\Delta|_{j+1}$ as a diary.
Extend  $e$  to $\Delta|_{j+1}$ as follows:
For $d\in \Delta_{j+1}$, 
let $e(d)$ equal  the $u\in D$ of length $|u_{j+1}|$ such that $u\re(|u_j|+1)=\lambda(d)$.
In particular, $e(d_{j+1})=u_{j+1}$.

If $j+1=n-1$, then we are done.
Otherwise, define  $\Delta_{j+2}$ 
as follows:
For each $u\ne u_{j+1}$ in $D$ of length $|u_{j+1}|$, $u$ has exactly one immediate successor in $D$ of length $|u_{j+2}|$.
For all $d\in \Delta_{j+1}\setminus\{d_{j+1}\}$,
define  $d^+=d^{\frown}k$, where $k\in 3$  satisfies ${e(d)}^{\frown}k\in\widehat{D}$.
Note that $u(|u_{j+1}|)=k$ for any/all $u\in D$ extending $e(d)$.
If $u_{j+1}$ is of type (a), then 
$d_{j+1}$ is of type (a) and we
 let $\Delta_{j+2}=\{{u_{j+1}}^{\frown}0  ,{u_{j+1}}^{\frown}2 \}\cup\{d^+:d\in \Delta_{j+1}\setminus\{d_{j+1}\}\}$.
If  $u_{j+1}$ is of type (b), then 
$d_{j+1}$ is of type (b)
and we
 let $\Delta_{j+2}=\{d^+:d\in \Delta_{j+1}\setminus\{d_{j+1}\}\}$.
If  $u_{j+1}$ is of type (c) (type (d), resp.), then 
$d_{j+1}$ is of type (c) (type (d), resp.),
and we
 let $\Delta_{j+2}=\{{u_{j+1}}^{\frown}1\}\cup
\{d^+:d\in \Delta_{j+1}\setminus\{d_{j+1}\}\}$.
This decides the type of critical node $d_{j+1}$ is in $\Delta$ and defines $\Delta_{j+2}$.
Note that the lexicographic order preserving map from $\Delta_{j+2}$ to $\{u\in D:|u|=|u_{j+2}|\}$ is a bijection.

At the end of $n$ many steps, we have constructed a diary $\Delta$ and a diary embedding $e:\Delta\ra D$.
At each stage in the construction, the nodes in $\Delta$ were uniquely determined by $D$.
In particular, $e[\Delta]=D$ which is the unique diary-shaped subset of  $A$ induced by $C$.
\end{proof}

Finally, we are equipped to prove the Main Theorem of the paper. 
\ref{thmmain}

\begin{thmMain}
For each finite chain $\mathbf{C}$ in $\Psi$, 
$\BRD(\mathbf{C},\Psi)<\infty$.
In particular, $\BRD(\mathbf{C},\Psi)$ is no more than the number of diaries for $\mathbf{C}$.
\end{thmMain}

\begin{proof}
Fix a finite chain $\mathbf{C}$, and suppose
$\chi$ is a coloring of the copies of $\mathbf{C}$ in $\Psi^*$.
There are finitely many diaries encoding copies of $\mathbf{C}$ in $\bS$; let $\mathcal{D}$ be the collection of these diaries.
For each  diary $\Delta\in\mathcal{D}$, $\chi$ induces a coloring of $\Emb(\Delta,\bS)$, which we again denote by $\chi$. 
Apply Theorem \ref{thm.FBRD} finitely many times,  once for each  diary in $\mathcal{D}$, to 
obtain a $\phi\in \Emb(\bS,\bS)$  so that 
 for each  $\Delta\in\mathcal{D}$, 
the set 
$\{\phi\circ g:g\in\Emb(\Delta,\bS)$ is monochromatic for $\chi$.

Let $S=\phi(\bS)$.
By  Lemma \ref{lem.existsaac}, there is an almost antichain $A\sse S$ encoding a copy of $\Psi^*$.
By 
Theorem \ref{thm.Asubsetsdiaries}, each  copy of $\mathbf{C}$ in $\Psi^*|_A$ 
is the image of some unique 
 diary embedding.
Thus, $\chi$ takes at most $|\mathcal{D}|$ many colors on the copies of $\mathbf{C}$ in $\Psi^*|_A$.
\end{proof}


\subsection{Chains of length two. }\label{subsec.ex7}

We now show that the big Ramsey degree for chains of length two in $\Psi$ is at most seven.

\begin{cor}\label{cor.BRD7}
Chains of length two in $\Psi$ have big Ramsey degree at most seven.
\end{cor}

\begin{proof}
Let $\mathbf{C}$ be a chain of length two in $\Psi$.
By Theorem \ref{thmmain},  $\BRD(\mathbf{C},\Psi)$ is at most the number of diaries for $\mathbf{C}$.
We utilize the construction in Theorem \ref{thm.Asubsetsdiaries} to show that there are seven such diaries.

Let $A\sse\bS$ be an almost antichain representing a copy of $\Psi^*$, and let 
 $C=\{a,b\}\sse A$ be  two coding nodes   representing a chain in $\Psi^*|_A$, where $|a|<|b|$.
Let $D$ be the diary-shaped subset of $A$ induced by $C$, as constructed in Theorem \ref{thm.Asubsetsdiaries}.
Note that $b$
 must be  a terminal coding node, and hence is of type (b).
The least node in $D$ must be  either a splitting node (so of type (a)) or a non-terminal coding node (so of type (c)).
Let $s=a\cap b$.
The only way $s$ can be a node of type (c) is Case A.\ below.
In all other cases, both $a$ and $b$ are terminal coding nodes and 
 $s$ is a splitting node in $D$.
We group the diary shapes into three categories. 

\begin{enumerate}
    \item[A.] 
When $a$ is not a terminal coding node, then $s=a$,  $a^{\frown}1\sse b$, and there is  only one possible diary.
\begin{enumerate}
\item[(1)]
The diary is 
  $\Delta=\{\langle\rangle, \langle 1\rangle\}$, where $\langle\rangle$ is a coding node of type (c) and $ \langle 1\rangle$ is a coding node of type (b).
\end{enumerate}
\end{enumerate}
\vskip.1in
From now on both $a$ and $b$ are terminal coding nodes an $s$ is a splitting node.
Hence, $\langle\rangle$ will be a splitting node (so of type (a)) in each of the following diaries.
\begin{enumerate}
   \item[B.]
When $\theta(a)=\theta(b)$, then 
there are  two possible diaries.
    \begin{enumerate}
   \item[(2)]
$\Delta=\{\langle \rangle,    \langle 0\rangle, \langle 2\rangle, \langle 0,0\rangle\}$, where 
$ \langle 2\rangle$ and $ \langle 0,0\rangle$ 
are both coding nodes of type (b). 

\item[(3)]
$\Delta=\{\langle\rangle, \langle 0\rangle, \langle 2\rangle,\langle 2,2\rangle\}$, where
$\langle 0\rangle$ and  $\langle 2,2\rangle$ are both coding nodes  of type (b).

\end{enumerate}
\vskip.1in
    \item[C.] 
When $\theta(a)\ne\theta(b)$, then  whichever of $a$ or $b$ is in the leftmost path in $D$  must extend  at least one node of type (d).
Let $d$ denote the maximum initial segment of the leftmost coding node in $D$ such that $\theta(d)=\theta(s)$.
There are four possible diaries.
\begin{enumerate}
\item[(4)]
If the leftmost coding node of $D$ is shorter than the rightmost, 
then  $D$ has critical nodes 
$\{s,d,a,b\}$, where
$s^{\frown}0\sse d\sse a$ and  $s^{\frown}2\sse b$.
Then  $\Delta=\{ \langle\rangle, \langle 0\rangle, \langle 2\rangle, \langle 0,1\rangle,  \langle 2,2\rangle, \langle 2,2,2\rangle\}$, where $ \langle 0\rangle$ is of type (d), and  $\langle 0,1\rangle$ and $\langle 2,2,2\rangle$ are of type (b).
\end{enumerate}

\noindent From now on, we assume the leftmost coding node of $D$  is longer than the rightmost.
\begin{enumerate}

\item[(5)]
If $|a|<|d|$, then 
$D$ has critical nodes 
$\{s,a,d,b\}$, where
$s^{\frown}0\sse d\sse b$ and  $s^{\frown}2\sse a$.
Then 
 $\Delta= \{ \langle\rangle, \langle 0\rangle, \langle 2\rangle, \langle 0,1\rangle, \langle 0,1,0\rangle\}$, 
where $ \langle 0,1\rangle$ is of type (d) and 
$\langle 2\rangle$ and $\langle 0,1,0\rangle$ are of type (b).

\item[(6)]
$|d|<|a|$ and 
$\theta(b)=\theta(b\re|a|)$,
then $D$ has critical nodes $\{s,d,a,b\}$ and
the diary is $\Delta =\{ \langle\rangle, \langle 0\rangle,  \langle 2\rangle, \langle 0,1\rangle,
\langle 2,2\rangle, \langle 0,1,0\rangle\}$, 
where $\langle 0\rangle$ is of type (d) and 
$\langle 2,2\rangle$ and $ \langle 0,1,0\rangle$ are of type (b).

\item[(7)]
$|d|<|a|$ and 
$\theta(b)\ne\theta(b\re|a|)$, then 
$D$ has critical nodes $\{s,d,a,d', b\}$, where $d'$ is the maximal initial segment of $b$ such that 
$\theta(d')=\theta(d'\re|a|)$.
Then the diary is 
$\Delta= 
\{ \langle\rangle, \langle 0\rangle,  \langle 2\rangle, \langle 0,1\rangle,  \langle 2,2\rangle, \langle 0,1,0\rangle, \langle 0,1,0,1\rangle\}$,
where $\langle 0\rangle$ and  $\langle 0,1,0\rangle$
are of type (d), and $ \langle 2,2\rangle$ and $ \langle 0,1,0,1\rangle$
 are of type (b).
\end{enumerate}
\end{enumerate}
\end{proof}

Combining Corollary  \ref{cor.BRD7} and the result in \cite{CEW_EUROCOMB25} that there is a coloring of the chains of length two in $\Psi$ into seven colors such that each color persists in every subcopy of $\Psi$, we obtain the following.

\begin{cor}
The big Ramsey degree for chains of length two in the two-branching countable ultrahomogeneous pseudotree is exactly seven.
\end{cor}

\subsection{Further directions}\label{subsec.furtherdir}

We close by mentioning that 
the work in this paper was done in a manner so as to set up for  our forthcoming work.
One line is to prove that the upper bounds for  finite chains in the pseudotree proved in Theorem \ref{thm.FBRD} are exact.
Another is to develop topological Ramsey spaces whose elements represent copies of the pseudotree,  giving infinite-dimensional Ramsey theorems which have the additional property of  directly recovering the exact big Ramsey degrees for finite chains.
Lastly, the work in this paper sets the stage for big Ramsey degree investigations on pseudotrees related to all generalized Wa\.{z}ewski dendrites.

\bibliographystyle{amsplain}
\bibliography{referencesPsi}

\end{document}